# HAMILTON CIRCUITS IN GRAPHS AND DIGRAPHS

**HOWARD KLEIMAN**

## I. INTRODUCTION.

Let h be an arbitrary n-cycle in $S_n$. Assume $h'$ is another n-cycle in $S_n$. Then if $\sigma = (h)^{-1}h'$, $\sigma$ is an element of $A_n$, the alternating group of $S_n$. From elementary group theory, every element of $S_n$ can be written as a product of (not necessarily disjoint) 3-cycles. Thus, $h' = h\sigma_1\sigma_2...\sigma_r$ where each $\sigma_i$ (i = 1,2,...,r) is a 3-cycle. Defining $h_i = h\sigma_1\sigma_2...\sigma_i$, we now stipulate that each $h_i$ is an n-cycle. Furthermore, let $G_m$ be a graph with vertex set

$$V = \{1,2,...,n\}$$

and edge set

$$E = \{e_i \mid i = 1,2,...,m\}.$$

Now suppose that h represents a pseudo-hamilton circuit in $G_m$, say H, where H is a hamilton circuit in

$$G_{m''} = G_m \cup \{ f_i \mid i = 1,2,...,s \}$$

with each $f_i$ an element of $K_n - G_m$, where $K_n$ is the complete graph on n vertices. Assume that the vertices of H occur going in a clock-wise manner. Call the edges, $f_i$, ***pseudo-arcs of H*** and the vertices of H from which pseudo-arcs emanate ***pseudo-arc vertices*** of H. All other vertices of H are ***arc vertices.*** Furthermore, let H' be a hamilton circuit in $G_m$ corresponding to the n-cycle, h', in $S_n$. Then the sequence h, $h_1$, $h_2$, ... culminates in $h_r$ = h'. Each $h_{i+1}$ is obtained from $h_i$ by multiplying it by the 3-cycle, $\sigma_i = (a_1\ a_2\ a_3)$. Applying $h_i$ to point $a_j$ yields

$$h_i(a_j) = a_{j+1} \text{ for } j = 1,2$$

while

$$h_i(a_3) = a_1$$

Thus, the arcs of $H_{i+1}$ are obtained from those of $H_i$ by replacing

$$(a_1, H_i(a_1)) \text{ by } (a_1, H_i(a_2)),$$
$$(a_2, H_i(a_2)) \text{ by } (a_2, H_i(a_3))$$
$$(a_3, H_i(a_3)) \text{ by } (a_3, H_i(a_1)).$$

We now place the following condition on the three new arcs:

***at most one of them is a pseudo-arc.***

**If our graph is random and n to inf, although we must sometimes backtrack when we're applying an algorithm to a directed graph, probabilistically, the total number of edges lying in G $_m$ increases in the sequence**

$$H_0, H_1, H_2, \ldots, H_i, \ldots .$$

**In this paper, we give algorithms for obtaining a hamilton circuit covering a number of cases. In particular, using Theorems ABKS, Frieze-ABKS and A Blocking Theorem, our paper answers two conjectures of Frieze in the affirmative:**

**(1) D $_{2-in,2-out}$ almost always contains a hamilton circuit as n to inf.**

**(2) R $_3$, the regular 3-out graph on n vertices almost always contains a hamilton circuit as n to inf.**

**Furthermore, we give necessary and sufficient conditions for an arbitrary graph to have a hamilton circuit. We then prove in Theorem H that if an arbitrary random graph of degree n contains a hamilton circuit, Algorithm G can almost always obtain a hamilton circuit in O(n $^{4.5}$ (logn) $^4$ ) running time with probability greater than**

$$\left(1 - \frac{1}{n^{3n^2}}\right)^2$$

**We then use Algorithm D to prove a corresponding theorem for random directed graphs. We also make the following conjecture:**

**CONJECTURE A. Let G* be an arbitrary graph that contains a hamilton circuit. Then Algorithm G always obtains a hamilton circuit in G* in polynomial time.**

## II. PRELIMINARIES.

**A graph, G, is *m-connected* if, given any pair of vertices u and v, there exist m disjoint paths from u to v. A graph which is 1-connected is *connected.* A directed graph which is at least 1-connected is *strongly- connected.* Let N = $\binom{n}{2}$. In [8], Erdös and Renyi gave the probability**

$$\frac{1}{\left(\frac{N}{k_n}\right)}$$

**for randomly choosing a graph, G, with n vertices and k $_n$ edges. They then gave the following limit distribution as n to inf, for the event E(n, k $_n$ ) " The random graph, G, containing n vertices and k $_n$ edges is 2-connected." where**

**k $_n$ = $\left(\frac{1}{2}\right)$(n)(logn + loglogn) + c $_n$ :**

$$\lim P(E(n,k_n)) = \begin{cases} 0 \text{ if } c_n \text{ to - inf}, \\ \exp(e^{-2c}) \text{ if } c_n \text{ to } c, \\ 1 \text{ if } c_n \text{ to inf.} \end{cases}$$

**It follows that if G is 2-connected, then each vertex is at least of degree 2. In [20], Palásti proved the following: Let D be a random directed graph, $D_{n,N}$, where each edge is equally likely to be chosen from among all edges in $KD_n$, the complete directed graph on all vertices (including all loops). Then if N = $N_c$ where $N_c$ = [nlogn + c] and c is an arbitrary number, the probability that $D_{n,N}$ is strongly connected has the probability**

$$\lim_{n \text{ to inf}} P_{n,N_c} = \exp(-2e^{-c})$$

**It follows that as c approaches a very large positive number, the above limit approaches 1. In [17], Kómlos and Szemerédi proved the following theorem:**

**THEOREM KS. Let the edges in $K_n$ be numbered $e_1$, $e_2$, ... , $e_N$ where $N = \binom{n}{2}$. Suppose that each edge in $K_n$ has been chosen in the following manner: The first edge has been chosen randomly with probability $\frac{1}{\binom{n}{2}}$, the second edge has been chosen randomly from the remaining edges with probability $\frac{1}{\binom{n}{2}-1}$, etc. We stop this process at the first instant when every valence is at least 2, say when m edges have been chosen. Then**

$$\lim_{n \text{ to inf}} \Pr(G_m \text{ is hamiltonian.}) = \lim_{n \text{ to inf}} \Pr(\delta(G_m) \geq 2) = 1$$

**(Theorem KS also appears to have been proven by Ajtai, Kömlos and Szeremédi in [1]).**
**In [3], Bollobás proved the following:**

**THEOREM ABKS. Let $\{e_1, e_2, ..., e_N\}$ (where $N = \binom{n}{2}$) be a random permutation of $K_n$.**

**If**
**$G_m = \{V, \{e_1, e_2, ..., e_N\}$, while**
**m" = {min( m: $\delta(G_m) \geq 2$}, then**

$$\lim_{n \text{ to inf}} \Pr(G_{m''} \text{ is hamiltonian}) = 1$$

**We call the type of random graph constructed by Bollabás a *Boll graph.* An analogous theorem for directed random graphs, Theorem Frieze-ABKS, was proven by Frieze in [12]. Henceforth, we assume that the set of vertices of each graph or directed graph is V. We define the randomness of our choice of edges as Boll does in [3]:**
**Let**

$$p(E) = \{e_i \mid i = 1,2,\ldots,\frac{n(n-1)}{2})$$

be a random permutation of the edges in $K_n$ . Then $G_m$ is the random graph with edges

$$E = \{e_i \mid i = 1,2,\ldots,m\}$$

where m < $\frac{n(n-1)}{2}$ . Define k(G) to be the *vertex connectivity* of G.

Let $\mathbb{N}$ be the set of natural numbers.

In [5], Bollobás proved the following theorem:

THEOREM K. Let Q be a monotonically increasing property of graphs and t a function defined by

$$t(Q) = t(Q;G_{m''}) = \min(m: G_m \text{ has } Q).$$

Then, given d ε $\mathbb{N}$, as n to inf,

$$t(\delta(G_{m''}) \geq d) = t(k(G_{m''}) \geq d) .$$

It follow that if a Boll graph, $G_{m''}$, has each vertex of degree at least 2, then $G_{m''}$ is almost always 2-connected. In [21], Wormald proved that almost all random graphs on n vertices which are of degree r are r-connected as n to inf. In [14], Frieze, Jerrum, Molloy and Wormald proved that almost all random regular graphs on n vertices which are of degree 3 have hamilton circuits as nto inf. Let $D_m$ be a random, directed in which m arcs have been randomly chosen to emanate from each vertex. If we change each arc into an unoriented edge, then the resulting graph, $R_m$, is a *regular m-out-degree graph.* Let $D_{i,o}$ be a directed graph on n vertices constructed in the following manner: Randomly choose i edges out of each vertex and o edges into each vertex. Then $D_{i,o}$ is an *i-in, o-out directed graph.* In [11], Fenner and Frieze proved that $R_m$ is m-connected and that $D_{i\text{-}in,\ o-out}$ is strongly connected. These results are necessary to apply Algorithms G and D, respectively, to $R_3$ and $D_{2-in,2-out}$ . Assume that $G_{m''}$ is a random graph satisfying the hypotheses of Theorem ABKS. Then Algorithm G yields a hamilton circuit with probability approaching 1 as n to inf. Similarly, let $D_{m''}$ be a random directed graph that satisfies the hypotheses of Theorem Frieze-ABKS. Then Algorithm D obtains a hamilton cycle in $D_{m''}$ with probability approaching 1 as n to inf. V and

$E = \{e_i \mid i = 1,2,\ldots,m\}$ define the respective vertices and arcs of a random directed graph, $D_p$, in which each arc is chosen with a fixed probability, p. In [2], Angluin and Valiant describe an O(nlogn) algorithm, A, such that

$$\lim_{n \text{ to inf}}(\Pr(A \text{ finds a hamilton circuit in } D\ .)) = 1 - n^{-\alpha}$$

where p = $\frac{c\log n}{n}$ with c dependent on α . Let

$$m = \frac{1}{2}(n\log n + n\log\log n + c(n))$$

where $G_m$ is a graph chosen randomly from among all graphs containing m edges. The latter was the definition of a random graph used by Erdös and Renyi in [8]. Bollobás,

**Fenner and Frieze used the same definition of a random graph in [4]. In [4], they gave an algorithm for obtaining a hamilton circuit in a random graph in $G_m$ with running time $O(n^{3+o(1)})$ .**

**In [19], McDiarmid proved that if $D_m$ is a random directed graph with**
**$m = n(\log n + c_n)$, then**

**lim (Pr($D_m$ is hamiltonian.)) to 1**
**n to inf**

**In [12], Frieze gave a sharp threshold algorithm with running time $O(n^{1.5})$ for obtaining a hamilton circuit as n to inf. In [13], Frieze and Luczak proved that as n to inf, $R_5$ almost always has a hamilton circuit. In [6], Cooper and Frieze proved that as n to inf, $D_{3-in,3-out}$ almost always has a hamilton circuit.**

**THEOREM A. Let G be a random graph (random directed graph) containing a pseudo-hamilton circuit (cycle) H. Then the probability that a pseudo-3-cycle is H-admissible is**

$$\frac{n-3}{2(n-2)}$$

**PROOF. Let H be a pseudo-hamilton circuit represented by equi-distantly-spaced points traversing the circle in a clock-wise manner. Now construct a random chord (1,H(j)) which represents an oriented edge or *arc* of G. Here**
**$1 < j \le n$. Thus, we cannot let j = 2, since (1,2) lies on the circle. Thus, the probability, $p_1$, that H(j) is chosen as the terminal point of the arc is $\frac{1}{n-2}$ . Next, randomly construct an arc, (j, H(k)). It is easily shown by constructing hσ that**

**S = {(1, H(j)), (j, H(k))}**

**defines an H-admissible pseudo-3-cycle, σ = (1 j k), if and only if the two arcs intersect in the circle H. The probability, p, that they intersect is**
**Pr(We randomly choose (1, H(j)).)**
**X Pr((j, H(k)) intersects H(j) in H.)**
**We are assuming that no arc chosen is an arc of the directed hamilton circuit H. Thus, if m = n-2 and j' = j-2,**

$$p_2 = \left(\frac{1}{n-2}\right)\left(\sum_{j=3}^{j=n}\frac{n-j}{n-2}\right) = \frac{(n-2)(n-3)}{2(n-2)^2} = \frac{n-3}{2(n-2)}$$

**The following theorem of W. Hoeffding is given in [15]:**

**THEOREM OF HOEFFDING. Let B(a,p) denote the binomial random variable with parameters a and p with**

**BS(b,c;a,p) = Pr(b ≤ B(a,p) ≤ c)**

**Then**

**(i) BS(0,(1 – α)ap; a, p) ≤ exp(- $\alpha^2$ ap/2)**

**(ii) BS((1 + α)ap, ∞; a, p) ≤ exp(- $\alpha^2$ ap/2)**

**where α is a small, positive number.**

**THEOREM B. Let v be a randomly chosen vertex of a Boll random graph, $G_m$. Then the probability that v has precisely two edges of $G_m$ incident to it is at most**

$$\frac{\log(cn(\log n)^2}{2n}$$

**PROOF. We first define hypergeometric probability.**
**Consider a collection of $N = N_1 + N_2$ similar objects, $N_1$ of them belonging to one of two dichotomous classes (say red chips), $N_2$ of them belonging to the second class (blue chips). A collection of r objects is selected from these N objects at random and without replacement. Given that**

$$X \varepsilon N, x \leq r, x \leq N_1, r-x \leq N_2,$$

**find the probability that exactly x of these r objects is chosen. Pr(X = x) is given by the formula**

$$\Pr(X = x) = \frac{\binom{N_1}{x}\binom{N_2}{r-x}}{\binom{N}{r}}$$

**where x objects are red and r – x objects are blue. Let v be a randomly chosen vertex of $G_m$. We wish to obtain the probability that exactly two edges in $G_m$ are incident to it. Let $N = N_1 + N_2$ where N is the sum of the degrees of all vertices in $K_n$ and $N_1$ is the degree of v in $K_n$. r equals twice the minimum number of edges in $G_m$ for which $G_m$ is 2-connected as n to inf. Thus,**

$$N = 2\binom{n}{2} = n(n-1)$$

**$N_1$ = the number of edges in $K_n$ incident to v**
**= n – 1**
**$N_2 = N - N_1 = n(n-1) - (n-1) = (n-1)^2$**
**r = the sum of the degrees of the vertices in $G_{m*}$**
**= at most [nlog(cn(logn)$^2$)]**
**x = 2**
**In r, we assume that c is a real number.**
**Then**

$$\Pr(X = 2) = \frac{\binom{n-1}{2}\binom{n^2-2n+1}{[n\log(cn(\ \log n)^2)]-2}}{\binom{n^2-n}{[n\log(cn(\log n)^2]}}$$

**From W. Feller [9], using the approximation of the hypergeometric distribution to the binomial distribution when N to inf, let p =** $\frac{n-1}{N} = \frac{1}{n}$ **yielding**

**Pr(X = 2) implies**

$$\text{B(2;N,p)} = \binom{[n\log(cn(\log n)^2)]}{2}\left(\frac{1}{n-1}\right)^2\left(1-\frac{1}{n-1}\right)^{[n\log(cn(\log n)^2]-2}$$

$$\approx [.5\log^2(cn(\log n)^2)]\exp(-\log(n\log(cn(\log n)^2))$$

$$\approx \frac{[\log(cn(\log n)^2]^2}{2(n\log(cn(\log n)^2)}$$

$$\approx \frac{\log(cn(\log n)^2)}{2n}$$

**concluding the proof.**

**COROLLARY TO THEOREM B. The probability that G contains more than** $(\log(cn(\log n)^2))^2$ **vertices of degree 2 approaches 0 as n to inf .**

**PROOF. Using Hoeffding's Theorem, let p =** $\frac{\log(cn(\log n)^2)}{2n}$ **while** **a =** $n\log(cn(\log n)^2)$**. Then from (ii) of his theorem, the probability, p'', that (1 + α)ap vertices are of degree 2 satisfies**

$$p'' < \exp\left(\frac{-\alpha^2(\log(\log(cn(\log n)^2)}{6}\right)$$

**which approaches 0 as n to inf. But for small α > 0,**

$$(1+\alpha)ap < (\log(cn(\log n)^2))^2,$$

**concluding the proof.**

**Define the *degree, v , of a directed graph* as** $\delta^-(v) + \delta^+(v)$**.**

**THEOREM C. Let** $D_{m''}$ **be a directed graph satisfying the hypotheses of Theorem Frieze-ABKS. Then, given a randomly chosen vertex, v, the probability that a unique arc emanates from v is no greater than** $\frac{\log n}{n}$ **as n to inf.**

**PROOF.** We again use hypergeometric probability. W.l.o.g., let $KD_n$ be the complete directed graph containing all arcs between any two vertices in V. Then let

N = the number of arcs in the complete directed graph on n vertices, $KD_n$
= n(n – 1)
$N_1 = N - N_1 = (n-1)(n-1) = (n-1)^2$,
r = the number of arcs in $D_{m''}$
= at most n(logn + c)
x = 1 .

Note. Since we are assuming that the hypotheses of the Frieze-ABKS Theorem are satisfied, if a vertex of $D_{m''}$ has precisely two arcs of $KD_n$ incident to it, then it must have precisely one arc terminating in it and one arc entering it. From hypergeometric probability,

$$\Pr(X = 1) = \frac{\binom{n-1}{1}\binom{(n-1)^2}{[n(\log n + c) - 1]}}{\binom{n^2 - n}{1}}$$

Again using the approximation of the hypergeometric distribution to the binomial distribution when N to inf, we obtain
Pr(X = 1) implies

$$B(1;N,p) = \binom{[n(\log n + c)]}{1}\left(\frac{1}{n-1}\right)\left(1 - \frac{1}{n-1}\right)^{[n(\log n + c)] - 1}$$
$$\to [\log n + \log c']\exp(-[\log n + \log c'])$$
$$\leq \frac{\log n}{n}$$

**COROLLARY TO THEOREM C.** The probability that there exist more than $2[(\log n)^2 + c\log n]$ arcs of $D_{m''}$ which are the unique arcs emanating from or terminating in a vertex approaches 0 as n to inf.

**PROOF.** The probability, p, that a randomly chosen vertex, v, of $D_{m''}$ has a unique arc emanating from it is at most $\frac{\log n}{n}$. The number of arcs in $D_{m''}$ is at most a = n(log n + c). The number of vertices is n. Thus, ap is at most $(\log n)^2 + c\log n$. From Hoeffding's Theorem (ii),

$$BS(1 + \alpha)ap, \infty; a,p) \leq \exp(-\frac{\alpha^2((\log n)^2 + c\log n)}{3}) \to 0$$

as n to inf . The same probability is true for the case where a unique arc terminates in v, concluding the proof.

**An *H-admissible product of two disjoint (pseudo) 2-cycles* (for short an *H-admissible POTDTC ),* say s = {(a,b), (c,d)}, occurs if and only if the vertices a, b, c, d traverse H in a clockwise manner in one of the following ways:**

**(i) a – c – b – d,**

**(ii) a – d – b – c .**

**It follows that if the vertices are placed on H in positions $\frac{2\pi j}{n}$ for j = 1,2, ... ,n , then, w.l.o.g., [a, H(b)] and [c, H(d)] are properly intersecting chords of H.**

**Before going further, we mention the following results given in Dickson [7]:**

**(i)** $\sum_{j=1}^{j=n} j = \frac{(n+1)(n)}{2}$

**(ii)** $\sum_{j=1}^{j=n} j^2 = \frac{n(2n+1)(n+1)}{6}$

**(iii)** $\sum_{j=1}^{j=n} j^3 = \left(\frac{n(n+1)}{2}\right)^2$

**THEOREM D. Let H be a pseudo-hamilton circuit of a random graph or a random directed graph G. Assume that G contains n vertices and that $e_1$ and $e_2$ are randomly chosen edges of G neither of which is an arc of H. Then the probability that $e_1$ and $e_2$ properly intersect (have no endpoints in common) is $\frac{n-3}{3(n-2)}$.**

**PROOF. W.l.o.g., let $e_1$ = (1,j). Consider the probability, p , that (1,j) properly intersects (r,s) where 2 ≤ r ≤ j-1 while j+1 ≤ s ≤ n. j can range over the domain [3,n]. It follows that given a specific value for j, the number of integral values of r is j-2: the edge, [r, H(r)], is not permissible by hypothesis; furthermore, the loop [r, r] is not an edge of $K_n$ . The number of possible values for s is n-j, namely, all vertices not contained in [1, j]. Thus, given a fixed value of j, the number of successes equals (n-j)(j-2). It follows that the total number of successes is**

$$\sum_{j=3}^{j=n}(n-j)(j-2)$$

**On the other hand, for a fixed value of j, the number of failures equals the number of possible values of r (j-2) multiplied by the number of possibilities for s. s cannot lie in the closed interval [n-j,n]. Furthermore, it must be distinct from r and H(r). Therefore, the number of possibilities for s is j-2. Therefore, the number of possibilties for failure is $(j\text{-}2)^2$. It follows that, using all values of j, the number of possibilities for failure is**

$$\sum_{j=3}^{j=n}(j-2)^2$$

**Thus, the probability of intersection is**

$$\frac{\sum_{j=3}^{j=n}(n-j)(j-2)}{\sum_{j=3}^{j=n}(n-j)(j-2)+(j-2)^2}$$

**Now let j' = j-2, m = n-2. Then the probability of success simplifies to**

$$\frac{\sum_{j'=1}^{j'=m} mj'-(j')^2}{m(\sum_{j'=1}^{j'=m} j')} = \frac{\frac{m^2(m+1)}{2}-\left(\frac{m(2m+1)(m+1)}{6}\right)}{\frac{m^2(m+1)}{2}} = \frac{m(m+1)(3m-2m-1)}{m(m+1)(3m)} = \frac{m-1}{3m}$$

$$= \frac{n-3}{3(n-2)}$$

**LEMMA 1. Let [a, H(b)] be an edge of a random graph, while H is a pseudo-hamilton associated with G. Let [b, H(c)] and [H(a), d] be edges in G neither of which contains an arc of H. Let m = n-2, j ' = j-2. Then if the minimum degree of G is 3, the following hold:**

**(a) The probability that *both* [b, H(c)] and [H(a), d] intersect [a, H(b)] is**

$$\frac{2m^2-3m+1}{6m^2}$$

**(b) The probability that [b, H(c)] *intersects* [a, H(b)], while [H(a), d] *doesn't intersect* [a, H(b)] is**

$$\frac{m^2-1}{6m^2}$$

**(c) The probability that [H(a), d] *intersects* [a, H(b)] , while [b, H(c)] *doesn't intersect* [a, H(b)] is**

$$\frac{m^2-1}{6m^2}$$

**(c) The probability that *neither* [H(a), d] *nor* [a, H(b)] intersects [a, H(b)] is**

$$\frac{2m^2+3m+1}{6m^2}$$

**PROOF. (a) Let a = 1 and H(b) be j. Since [1, j] doesn't lie on the directed circuit, H , the possible choices for j are 3,4, ... ,n . The possible choices for b are 2,3, ... , j-1. Thus, there are j-2 choices for both j and b. In order for [b. H(c)] to properly intersect [a, H(b)], H(c) must be one of the numbers j+1, j+2, ... ,n ,implying that there are n-j possibilities for a choice of H(c). On the other hand, H(a) = H(1) = 2. Given a fixed value for j, the probability that [H(a), d] intersects [a, H(b)] implies that d is one of the numbers j+1, j+2, ... ,n. Thus, there are n-j possibilities for d. Therefore,**

$$p = \left(\frac{1}{n-2}\right)\left(\sum_{j=3}^{j=n}\frac{(n-j)^2}{(n-2)^2}\right)$$

$$= \sum_{j'=1}^{j'=m}\frac{(m-j')^2}{m^3}$$

$$= \frac{m^3 - m(m+1)(m) + \dfrac{m(2m+1)(m+1)}{6}}{(m)^3}$$

$$= \frac{6m^3 - 6m^3 - 6m^2 + 2m^3 + 3m^2 + m}{6m^3}$$

$$= \frac{2m^2 - 3m + 1}{6m^2}$$

**which approaches $\frac{1}{3}$ as n to inf.**

**PROOF OF (b). In this case, given b, we have n-j choices for H(c). On the other hand, since [H(a), d] *doesn't* intersect [1,j], we have j-2 choices for d. It follows that the probability is**

$$\sum_{j=3}^{j=n-}\frac{(n-j)(j-2)}{(n-2)^3} = \sum_{j'=1}^{j'=m-}\frac{(m-j')(j')}{m^3}$$

$$= \sum_{j'=1}^{j'=m}\frac{m\,j' - (j')^2}{m^3}$$

$$= \frac{\dfrac{m(m+1)(m)}{2} - \dfrac{m(2m+1)(m)}{6}}{m^3}$$

$$= \frac{\mathbf{3}m^2 + \mathbf{3}m - \mathbf{2}m^2 - \mathbf{3}m - \mathbf{1}}{6m^2}$$

$$= \frac{m^2 - 1}{6m^2}$$

**PROOF OF (c). The proof is the same as that of (b).**

**PROOF OF (d). Since *neither* [H(a), d] *nor* [b, H(c)] intersect [a, H(b)], there are j-2 possibilities for d and j-2 possibilities for H(c). It follows that the probability for two successes and no failures is**

$$\sum_{j=3}^{j=n} \frac{(j-2)^2}{(n-3)^3} = \sum_{j'=1}^{j'=m} \frac{(j')^2}{m^3} = \frac{(2m+1)(m+1)}{6m^2} = \frac{2m^2 + 3m + 1}{6m^2}$$

**LEMMA 2. Let [a, H(b)] be defined as in LEMMA 1. Let [b, H(c)], [b', H(c')], [H(a), d], [H(a), d'] be edges belonging to a random graph G. Then the following hold:**

**(a) The probability that [b, H(c)], [H(a), d] and [H(a), d'] all intersect [a, H(b)] is**

$$p_3 = \frac{m^2 - 2m + 1}{4m^2}$$

**(d) The probability that precisely two of these edges intersect [a, H(b)] is**

$$p_2 = \frac{m^2 - 1}{4m^2}$$

**(e) The probability that precisely one of these edges intersects [a, H(b)] is**

$$p_1 = \frac{m^2 - 1}{4m^2}$$

**(f) The probability that none of the three edges intersects [a, H(b)] is**

$$p_0 = \frac{m^2 + 2m + 1}{4m^2}$$

**PROOF OF (a). Let m = n-2 and j' = j-2. Then**

$$p_3 = \sum_{j=3}^{j=n} \frac{(n-j)^3}{(n-2)^4}$$

$$= \sum_{j'=1}^{j'=m} \frac{(m-j')^3}{m^4}$$

$$= \sum_{j'=1}^{j'=m} \frac{m^3 - 3j'm^2 + 3m(j')^2 - (j')^3}{m^4}$$

$$= \frac{m^4 - \frac{3m^2(m+1)(m)}{2} + \frac{3m^2(m)(2m+1)(m+1)}{2} - \frac{(m+1)m^2}{4}}{m^4}$$

$$= \frac{12m^4 - 18m^4 - 18m^3 + 12m^4 + 18m^3 + 6m^2 - 3m^4 - 6m^3 - 3m^2}{12m^4}$$

$$= \frac{m^2 - 2m + 1}{4m^2}$$

**PROOF OF (b).**

$$p_b = \sum_{j=3}^{j=n} \frac{(n-j)^2(j-2)}{(n-2)^4}$$

$$= \sum_{j'=1}^{j'=m} \frac{(m-j')^2 j'}{m^4}$$

$$= \sum_{j'=1}^{j'=m} \frac{m^2 j' - 2m(j')^2 + (j')^3}{m^4}$$

$$= \frac{\frac{m^2(m+1)(m)}{2} - \frac{2m(m)(2m+1)(m+1)}{6} + \frac{(m+1)^2 m^2}{4}}{m^4}$$

$$= \frac{m^2 - 1}{12m^2}$$

**It follows that since there are three possibilities here, namely,**

**T[b, H(c)]F[H(a), d]F[H(a), d'],**
**F[b, H(c)]T[H(a), d]F[H(a), d'],**

**F[b, H(c)]F[H(a), d]T[H(a), d '],**

**each with the same probability of success,**

$$p_2 = \frac{m^2 - 1}{4m^2}$$

**PROOF OF (c).**

$$p_c = \sum_{j=3}^{n=n} \frac{(n-j)(j-2)^2}{(n-2)^4}$$

$$= \sum_{j'=1}^{j'=m} \frac{(m-j')(j')^2}{m^4}$$

$$= \sum_{j'=1}^{j'=m} \frac{m(j')^2 - (j')^2}{m^4}$$

$$= \frac{\frac{m(m)(2m+1)(m)}{6} - \frac{m^2(m+1)^2}{4}}{m^4}$$

$$= \frac{m^2 - 1}{12m^2}$$

**implying that**

$$p_1 = \frac{m^2 - 1}{4m^2}$$

**PROOF OF (d).**

$$p_0 = \sum_{j=3}^{j=n} \frac{(j-2)^3}{(n-2)^4}$$

$$= \sum_{j'=1}^{j'=m} \frac{(j')^3}{m^4}$$

$$= \frac{\frac{(m+1)^2(m)^2}{4}}{m^4}$$

$$= \frac{m^2 + 2m + 1}{4m^2}$$

**THEOREM E. Let a be a pseudo-arc vertex of a pseudo-hamilton circuit, H, of a random graph, $G_m$, where δ(G) ≥ 3. Then the probability that at least two H-admissible permutations containing a can be constructed is at least**

$$p = \frac{13m^4 - 8m^3 - 6m^2 + 1}{16m^4}$$

**where m = n – 2.**

**PROOF. Let $p_i$ (i= 0,1,2,3} be the probability that the edge [a, H(b)] is intersected by i edges of the set E = {[b, H(c)], [H(a), d], [H(a), d']} . We now consider the probability tree of all possibilities of intersections of [a, H(b)] with subsets of E. Before continuing, let**

**A = [a, H(b)], B = [b, H(c)], D = [H(a), d],**
**D' = [H(a), d']**

**Let X ∈ {B, D, or D'} . The statement, S , "X intersects A" is denoted by T(AX), while ~S is denoted by F(AX).**

**0 edges.**

**Only one branch of the tree satisfies this condition: F(AB)F(AD)F(AD'). From Lemma 2, (d),**

$$p_0 = Pr(F(AB)F(AD)F(AD'))$$
$$= \frac{m^2 + 2m + 1}{4m^2}$$

**(1 edge) 3 branches occur on the tree:**

**T(AB)F(AD)F(AD'), F(AB)T(AD)F(AD'), F(AB)F(AD)F(AD').**

**From Lemma 2, (c), the probability of each possibility is**

$$\frac{m^2 - 1}{12m^2}$$

**implying that**

$$\text{Pr(Precisely one edge intersects A.)} = \frac{m^2 - 1}{4m^2}$$

**(2 edges). 3 branches occur on the tree:**

**T(AB)T(AD)F(AD'), T(AB)F(AD)T(AD'), F(AB)T(AD)T(AD').**

**From Lemma 2, (b) ,the probability of each possibility is**

$$\frac{m^2 - 1}{12m^2}$$

**Thus,**

$$p_2 = \frac{m^2 - 1}{4m^2}$$

**(3 edges) One branch occurs on the tree:**
**T(AB)T(AD)T(AD').**

**From Lemma 2, (a),**

$$p_3 = \frac{m^2 - 2m + 1}{4m^2}$$

**Assume now that A' = [a, H(b')], B' = [b', H(c')],**
**D = [H(a), d], D' = [H(a), d']. If S' = {B', D, D'}, we obtain corresponding probabilities, $q_i$ (i = 0,1,2,3}, for the intersection of A' with the subsets of S', namely, $q_0$, $q_1$, $q_2$, $q_3$. The simplest way to obtain the probability, p, that at least two intersections occur with either A or A' or both A and A' is to obtain**
**p' = 1 – p. The only possible branches of the probability trees with respect to A and A' which yield at most one intersection are the following:**

**(i) 1 intersection of A, 0 intersections of A'.**
**(ii) 0 intersections of A, 1 intersection of A'.**
**(iii) 0 intersections of A, 0 intersections of A'.**

**Thus, $p' = p_0 q_0 + p_0 q_1 + p_1 q_0$, yielding**

$$p' = \frac{m^4 + 4m^3 + 6m^2 + 4m + 1}{16m^4} + \frac{2m^4 + 4m^3 - 4m - 2}{16m^4} = \frac{3m^4 + 8m^3 + 6m^2 - 1}{16m^4}$$

**It follows that**

$$p = \frac{13m^4 - 8m^3 - 6m^2 + 1}{16m^4} .$$

**LEMMA 3. Let D be a random directed graph with $\delta^+(D) \geq 2$, $\delta^-(D) \geq 2$. Let H be a pseudo-hamilton cycle of D. Then the probability of obtaining at least two H-admissible permutations containing an arbitrary pseudo-arc vertex, v, is the same as in the previous case.**

**PROOF. All we have to do is assume that v has two arcs emanating from it and H(a) has at least two arcs terminating in it.**

**Before beginning a sketch of the algorithm, we discuss the probability of obtaining an $H_i$-admissible pseudo-3-cycle in $R_3$, a regular 3-out graph obtained from the directed graph, $D_3$. Suppose that a is a pseudo-arc vertex of the pseudo-hamilton circuit, $H_i$. Let $H_{i''}$ be the graph consisting of $H_i$ together with all arcs of E symmetric to arcs of $H_i$. Suppose we randomly chose an edge incident to a, say [a, $H_i$(b)], lying in $R_3$ - $H_{i''}$, and then randomly choose an edge incident to b, say [b, $H_i$(c)], in $R_3$ - $H_{i''}$. The question then arises: May we assume that these two edges are actually chosen independently of each other? The answer is yes: Each of these edges was obtained from randomly chosen arcs of $D_3$. It is possible that the arcs in $D_3$ might be (a, $H_i$(b)) and (H(c), b). The crucial thing is that the probability that the two arcs intersect approaches $\frac{1}{2}$ as n to inf. But the corresponding edges form an $H_i$-admissible permutation if and only if they intersect. Thus, we may assume that randomly chosen edges of the form [a, $H_i$(b)], [b, $H_i$(c)] have, as n to inf, a probability of $\frac{1}{2}$ of defining an admissible 3-cycle. Similarly, edges of the form [a, $H_i$(b)], [c, $H_i$(d)] have a probability of $\frac{1}{3}$ of defining an admissible POTDTC. We may thus randomly choose edges from $R_3$ when applying ALGORITHM G to it. Definitions given in the introduction apply to the remainder of the paper. Propositions 3.1 – 3.4 in Angluin and Valiant [2] prove that we may equivalently define a random graph (random directed graph) either in terms of**

**(i)** **a random choice of the graph (directed graph) from among all graphs with m edges (arcs),**

**or**

**(ii)** **from among the set of all graphs (directed graphs) on n vertices where each edge (arc) is chosen with probability** $\frac{m}{\binom{n}{2}}$ ($\frac{m}{2\binom{n}{2}}$).

**Thus, these definitions may be used interchangeably. If (a, b) is a directed arc in a directed graph, then (b, a) is an arc *symmetric* to (a, b). A random edge chosen in Algorithm G of the next section is always assumed to be incident to a fixed vertex, say v. In that sense, *it behaves like an arc emanating from v.* Let**

$$h = (a_1\ a_2\ \dots\ a_n)$$

**denote an n-cycle of $S_n$. Suppose we apply an H-admissible 3-cycle, s, to h to obtain h' = hs where s = (a b c). Then**

$$h' = (ah(b) \dots ch(a) \dots bh(c) \dots )$$

**Call this *abbreviated representation* of h'**

$$a^* = (ah(b) \dots ch(a) \dots bh(c) \dots)$$

**since we have omitted only subpaths belonging to h. We define a* to be an *abbreviation of h'*. Now let H be the pseudo-hamilton circuit corresponding to the n-cycle h in $S_n$, while σ is the pseudo-3-cycle corresponding to s. Then the term H-admissible is used** interchangeably with h-admissible. In particular, the pseudo-hamilton circuit, H, corresponding to h and representable by

A = [a,H(a), ... , b,H(b), ... , cH(b), ... ]

is replaced by the ***abbreviation***

A' = [a,H(b), ... ,cH(a), ... ,bH(c), ... ]

representing H' . As an example, if

h = (1 2 3 4 5 6 7 8 9 10 11 12), s = (1 4 8)

then h' = hs = (1 5 6 7 8 2 3 4 9 10 11 12)
= (1h(4) ... 8h(1) ... 4h(8) ... )

Thus, h' can be represented by the abbreviation

a' = (1h(4) ... 8h(1) ... 4h(8) ... )

Correspondingly, the pseudo-hamilton circuit, H', is completely determined by the abbreviation

A' = (1H(4) ... 8H(1) ... 4H(8) ... )

since the remaining points of H' all occur in the order in which they occurred in H. In particular, s maps aH(a) into aH(c), bH(b) into bH(c), and cH(c) into cH(a) where a = 1, b = 4, c = 8. Essentially, we are partitioning H into three subpaths which are joined together to form the pseudo-hamilton circuit H'. In this example, s is the ***permutation associated with the abbreviation A .*** Now let s" = (2 6)(3 7) be applied to h to obtain

h" = hs" = (1 2 7 8 5 6 3 4 9 10 11 12)
= (1 2 h(6) 8 h(4) 6 h(2) 4 H(8) ...)

Here we are partitioning h into four subpaths joined together to form h". In general, the format for an abbreviation using an h-admissible POTDT,
s" = (a c)(b d) , is

hs" = (a h(c) ... d h(b) ... c h(a) ... b h(d) ...)

Before going on, we define a ***rotation.*** Let

$$S = [v_i, v_{i+1}, \dots, v_{i+k}]$$

be a subpath of a pseudo-hamilton circuit, H , where $v_i$ is a pseudo-arc vertex of H. Assume that $[v_i, v_{i+k}]$ is an edge lying in G – H. Then

$$S' = [v_i, v_{i+k}, v_{i+k-1}, \dots, v_{i+1}]$$

is a subpath of a pseudo-hamilton circuit, H', where S' replaces

$$S = [v_i, v_{i+1}, \dots, v_{i+k}]$$

of H to yield H'. This procedure is called a ***rotation with respect to*** $v_i$ ***and*** $v_{i+k}$. Suppose we represent a subpath, S, of H containing in its interior only vertices of degree 2, by a new vertex (not in G) identified by its first vertex, $v_\alpha$ and last vertex, $v_\beta$. This new vertex, denoted by $v_\alpha v_\beta$, is called a ***2-vertex.*** In a rotation, the order of the vertices of a 2-vertex is changed. Thus, if a 2-vertex is $v_\alpha v_\beta$ ***before*** a rotation containing it, it becomes $v_\beta v_\alpha$ after the rotation. Let the new hamilton circuit obtained from replacing subpaths of the form S by 2-vertices be H'. Rotations were used by Bollobás, Fenner and frieze in [4]. We note that, by construction, the number of pseudo-arc vertices of H' is never greater than the number in H. H-admissibility of a 3-cycle requires that [a, c] and [b, d] properly intersect in a circle along which the vertices of H are equally spaced. On the other hand, H-admissibility of a POTDT, (a c)(b d), requires that [a, c] and [b, d] properly intersect. We will search in depth in each iteration on up to $(\log n)^2$ permutations. Thus, using abbreviations when testing for H-admissibility, we need only consider *at most* $24(\log n)^3$ points for the first iteration, $48(\log n)^3$ points for the second one, ... , 48rlogn points for the r-th one. Thus, in r iterations, we randomly go through $24(\log n)^3$ points. Suppose a recalculation of $H_i$ occurs between every $n^{\frac{1}{2}}$ iterations. Then (as will be shown in more detail in IV), it requires $O(n(\log n)^3)$ running time to go through $n^{\frac{1}{2}}$ iterations. On the other hand, it requires O(i) running time (r. t.) to construct a rotation in an abbreviation containing i points. It follows that it requires at most O(nlog n) running time to use a rotation in each of $n^{\frac{1}{2}}$ iterations. It follows that in O(nlog n) iterations, it requires at most $O(n^{1.5}(\log n)^4)$ r. t. to complete these operations. The latter will also be the r. t. for both Algorithm G and Algorithm D in the next section. The algorithm essentially consists of sequentially obtaining a sequence

$$H_0, A_1, \ldots A_{[n^5]}, H_{[n^5]}, A_{[n^5]+1}, \ldots, A_{2[n^5]}, H_{2[n^5]}, \ldots$$

in which, using Algorithm G for graphs, the number of pseudo-arc vertices is a monotonically decreasing function. Using Algorithm D for directed graphs, we are required to backtrack from some iterations. However, using a large number of iterations, the number of successful iterations becomes considerably greater than the number of failures. Thus, we here also replace pseudo-arc vertices by arc vertices. In general, as n to inf, we approach a hamilton circuit in the former case, and a hamilton cycle in the latter case. Definitions and examples of the following data structures comes from Knuth [16]: Given m entries, using a ***balanced, binary search tree***, we can, respectively, locate, insert, or delete any element, or rebalance the tree in O(logn) r. t.. In this paper, a ***LIFO, double-ended queue*** – henceforth called a ***queue*** – is a linear list in which all insertions are made at the beginning of the list, while deletions may be made at either end of the list. Finally, in order to improve our algorithm when searching for a hamilton circuit in $G_{m''}$ or $D_{m''}$, we make the following definitions: Let

$$S_\alpha = [v_\alpha, v_{\alpha_1}, v_{\alpha_2}, \ldots, v_{\alpha_i}, v'_\alpha]$$

be a subpath of $G_{m''}$ where each vertex, $v_{\alpha_j}$ (j = 1,2, ... ,i) has precisely two edges incident to it, while $v_\alpha$ and $v'_\alpha$ each has at least three edges incident to it. For each $S_\alpha$ in $G_{m''}$, define $v_\alpha v'_\alpha$ as a ***2-vertex*** of a graph $G'_{m''}$. Here $v_\alpha v'_\alpha$ replaces the subpath $S_\alpha$ and vertices contained in $S_\alpha$. In constructing $G'_{m''}$, we delete any edge of $G_{m''}$ which connects $v_\alpha$ and $v'_\alpha$. Our reason for doing this is that such an edge cannot lie on any hamilton circuit of $G_{m''}$ since the set
$\{v'_\alpha, v_\alpha, v_{\alpha_1}, v_{\alpha_2}, \ldots, v_{\alpha_i}\}$ forms a cycle in $G_{m''}$. Furthermore, we may assume that if there exist edges $[v'_\alpha, v]$, $[v, v_\alpha]$ in $G_{m''}$, then only one of the edges can lie on a hamilton circuit for essentially the same reason as in the previous case. In fact, for very large n, we can prove that it almost always never occurs that $v_\alpha$ and $v'_\alpha$ are both incident to the same vertex in $G_{m''}$. Our reasoning is a follows: By construction, $G_{m''}$ has fewer edges than $G_m$. But there are two equivalent ways of defining $G_m$:

(1) We can randomly choose $G_m$ from among all graphs containing m edges.

(2) We can choose each edge with probability

$$p = \frac{m}{\binom{n}{2}}.$$

Since m = $\frac{1}{2}$n(logn + loglogn + c) where c doesn't approach - inf, for very large n, p < $\frac{(\log n)^{1+\varepsilon}}{n-1}$ = p' where $\in > 0$. Thus, given a fixed subpath, $S_\alpha$, if $[v, v_\alpha]$ lies in $G_{m^*}$, the probability that $e_1 = [v, v'_\alpha]$ also lies in $G_{m^*}$ is less than p'. Call E the event " There exists no vertex, v, such that two edges incident to v are incident to $v_\alpha$ and $v'_\alpha$, respectively." It follows that the probability that $e_1$ ***doesn't*** lie in $G_{m''}$ is at least $1 - p'$. From Hoeffding's Theorem, the probability that more than

$$\frac{n(\log n + \log\log n + c)(\log n)^{1+\varepsilon}}{n-1} < 2(\log n)^{2+\varepsilon}$$

edges of $G_{m''}$ are incident to a fixed vertex approaches 0 as n to inf. Thus, the probability that more than $4(\text{logn})^{2+\alpha}$ edges of $G_{m''}$ are incident to either $v'_\alpha$ or $v_\alpha$ approaches 0 as n to inf. It follows that the probability ~E is at least

$$(1 - p')^{4(\log n)^{2+\varepsilon}}$$

From the Corollary to Theorem B, the probability that there exist more than $M = (\text{logcn}(\text{logn})^2)^2$ vertices of degree 2 in $G_{m''}$ approaches 0 as n to inf. Thus, the number of possible subpaths, $S_\alpha$, is no greater than M. It follows that the probability that ***no*** vertex, v, in $G_{m''}$ has the property that $[v, v_\alpha]$ and $[v, v'_\alpha]$ ***both*** belong to $G_{m''}$ is at least

$$(1 - \frac{(\log n)^{1+\varepsilon}}{n-1})^{4(\log n)^{2+\varepsilon}(\log(cn(\log n)^2)^2)}$$

$$\approx \exp(-\frac{(\log n)^{5+\varepsilon}}{n-1})$$

$$\to 1$$

as n to inf. In any event, even if a vertex, v, exists in an arbitrary graph such that event E is true, both $[v, v_\alpha]$ and $[v, v'_\alpha]$ can't both occur on the same hamilton circuit, $H'_i$: If e = $[v_\alpha v'_\alpha, v]$ lies on $H'_i$ where w.l.o.g. an edge incident to $v_\alpha$ precedes e on $H'_i$, we must choose an edge incident to v in $G'_{m''}$, say e', which doesn't lie on $H'_i$, precluding a choice of e. It follows that e' = [v, v'] where v' is either a vertex of V other than $v_\alpha$ or $v'_\alpha$, or else it is a 2-vertex containing neither $v_\alpha$ nor $v'_\alpha$. In general, when choosing an edge incident to a 2-vertex, all we need do is to determine whether an edge exists which is terminating in or emanating from the 2-vertex. For instance, if $v_\alpha v_\beta$ has an edge of $H'_i$ incident to $v_\alpha$, then we must choose an edge incident to $v_\beta$. On the other hand, if *no* edge incident to either vertex lies on $H'_i$, then we can choose an edge incident to *either* vertex. However, if our 2-vertex is written, $v_\alpha v_\beta$, and we choose an edge incident to $v_\beta$, we must rewrite our 2-vertex as $-v_\alpha v_\beta$. This indicates that the subpath represented by the 2-vertex goes from $v_\beta$ to $v_\alpha$. One other type of problem may occur. Suppose that two subpaths, say $S_\alpha$ and $S_\eta$, have an endpoint in common, say $v_{\alpha'} = v_\eta$. In this case, we delete all edges incident to $v'_\alpha$ except $[v_{\alpha_i}, v'_\alpha]$ and $[v'_\alpha . v_\eta, v_{\eta_1}]$. It follows that we obtain one, larger subpath, say $S_{\alpha''}$, of the form

$$[v_\alpha ,v_{\alpha_1} ,...,v_{\alpha_i} ,v'_\alpha = v_\eta ,v_{\eta_1} ,...,v_{\eta_j} ,v'_\eta]$$

Now consider the case when we're working on a directed graph. W.l.o.g., assume that it is $D_{m''}$. Let

$$S_\beta = \{v_\beta, v_{\beta_1}, v_{\beta_2}, ... ,v_{\beta_j} ,v'_\beta]$$

be a directed subpath of $D_{m''}$ whose interior vertices are all of degree 2 (in the sense that precisely one arc enters and one arc leaves each of these vertices). Now suppose that there exists a unique arc, a = $(v_\beta, v'_\beta)$, which is the unique arc terminating in $v'_\beta$. But then we can delete all arcs emanating from $v_\beta$ other than a. Similarly, if a is a unique arc emanating from $v_\beta$, no other arc terminating in $v'_\beta$ can lie on a hamilton circuit of $D_{m''}$. Thus, we can delete all arcs terminating in $v'_\beta$. In the case of a, our "subpath", $S_\beta$, consists of a single arc. We now note some differences from the unoriented case. First, $S_\beta$ is a *directed* subpath. Thus, its initial vertex *always* is $v_\beta$; its last vertex is *always* $v'_\beta$. As in the previous case, neither the arc $(v'_\beta, v_\beta)$ nor the arc $(v_\beta, v'_\beta)$ can lie on any hamilton cycle of $D_{m''}$. Thus, if they exist in $D_{m''}$, we can delete them. Using this information we can replace the subpaths, $S_\beta$, by 2-vertices of the form $v_\beta v'_\beta$ to form the contracted, directed graph, $D'_{m''}$. We don't have to concern ourselves with the arcs $(v, v_\beta)$ and $(v'_\beta, v)$ since they remain distinct arcs in $D'_{m''}$. As we shall see in the next section, in both Algorithm G and Algorithm D, we construct abbreviations each of which represents $[n^{.5}]$ iterations.

(ii) If an $H_i$-admissible permutation yields a new pseudo-arc vertex, say $a_1$ ,then the iteration contains a rotation starting at the pseudo-arc vertex, say from

$$S = \{a_1, a_2, ... , a_i, ... \}$$

to

$$S' = \{a_1, a_i, a_{i-1}, ... , a_2, a_{i+1}, a_{i+2}, ... \}$$

(iii) We change the signs of all vertices in the rotation except the first one. This indicates that we are traversing $H_i$ in a counter-clockwise manner. Thus, $S'$ should be written

$$\{a_1, -a_i, -a_{i-1}, -a_{i-2}, ... , -a_2, a_{i+1}, a_{i+2}, ... \}$$

In doing so, we must keep in mind that if $a_j$ is a pseudo-arc vertex in S, then $[a_j, a_{j+1}]$ is a pseudo-arc. After the rotation, $[a_{j+1}, a_j]$ is still a pseudo-arc. Thus, $a_{j+1}$ becomes a pseudo-arc in $S'$, a subpath of $H_{i+1}$. On the other hand, $[a_j, a_{j-1}]$ is a pseudo-arc only if $[a_{j-1}, a_j]$ was one in S. Thus, $a_j$ is not necessarily a pseudo-arc of $H_{i+1}$. When we have obtained a hamilton circuit, $H'$, in a contracted graph, $G'_{m''}$ or $D'_{m''}$, we substitute the subpaths defined by the 2-vertices into the respective graph (directed graph) to obtain a hamilton circuit (hamilton cycle) in the original graph (directed graph).

**A BLOCKING THEOREM.** The following hold:

(i) As n to inf, $R_3$ contains no blocking subgraph which would make it impossible to randomly go through each vertex in V.

(ii) As n to inf, $D_{2-in,2-out}$ contains no blocking subgraph which would make it impossible to randomly go through each vertex in V.

**PROOF.** (i) In [13], Frieze and Luczak proved that as n to inf, with probability approaching 1, $R_5$ contains a hamilton circuit. Therefore, $R_3$, contains fewer arcs randomly chosen out of each vertex, v, of $D_3$ than $R_5$ has arcs randomly chosen out of v. Thus, $R_3$ contains no blocking subgraph as n to inf.

(ii) In [6], Cooper and Frieze proved that $D_{3-in,\,3-out}$ contains a hamilton circuit as n to inf. Therefore, $D_{2-in,2-out}$ contains no blocking directed subgraph as n to inf.

Note. $R_3$ is 3-connected, while $D_{2-in,\ 2-out}$ is strongly-connected.

## III. SKETCH OF ALGORITM

In general, if G is an unoriented graph, it is represented as a balanced, binary search tree whose first branches are numbered 1 through n together with respective counters which register the number of edges incident to each vertex. If G is a directed graph, it is represented by a balanced, binary search tree containing 2n primary branches where n is the number of vertices in $D_{m''}$: the first n branches represent arcs emanating from the respective vertices 1 through n; the second set of n vertices represents terminating in the respective vertices 1 through n. In order to cover all cases, $\{G_{m''}, D_{m''}, R_3, D_{2,2}\}$, call all graphs and directed graphs mentioned here, G. We generally use the word *edge* when discussing both edges and arcs. Since the only edges employed in the algorithm are those incident to a fixed vertex, they are essentially used as arcs. In the algorithm, starting with a randomly chosen initial pseudo-hamilton circuit, $H_0$ ,we successively construct new ones using $H_i$-admissible permutations, $s_i$ (i=0,1,2,3, ... ), to respectively obtain $H_i$ (i=1,2,3, ... ). In constructing $H_0$, if $\delta(G) > 2$, we must include all subpaths containing vertices of degree 2. If G is a graph, we randomly orient each subpath before placing it in $H_0$. The procedure for obtaining subpaths is explained in the construction of TWOPATHS further on. ***Since $G_{m''}$ satisfies the hypotheses of Theorem ABKS, from Bollobás [5], it is 2-connected. Similarly, since $D_{m''}$ satisfies the hypotheses of Theorem Frieze-ABKS, using Palásti [20], we intimate*** ***that*** $\boldsymbol{D_{m''}}$ ***is strongly-connected.*** **Thus, using Theorem E with respect to $G_{m''}$ and $D_{m''}$, we replace $G_{m''}$ by its contracted graph, $G_{m''}$ ,and replace $D_{m''}$ by $D_{m''}'$. Preliminary to constructing a pseudo-hamilton circuit in $G_{m''}'$ ($D_{m''}'$), we place each subpath, $S_\alpha$ ($S_\beta$), in a balanced, binary search tree called TWOPATH. If G = $G_{m''}$, we count the number of edges incident to each vertex. If a vertex has two edges incident to it, we place the vertex together with the edges incident to it in TWOPATH. After we have finished creating TWOPATH, we check to see if any end vertex, v', of the edge, $S_\alpha$ = [v.'s, v'], in TWOPATH is itself a vertex of degree 2. If it is, we extend $S_\alpha$ to form $S_{\alpha'}$ = {v.'s, v', v"}. In this manner, we construct the largest subpaths possible. We place each subpath constructed in TWOPATH. We follow an analogous procedure if G = $D_{m''}$. Once we have constructed subpaths, we construct G' from $G_{m''}'$ by replacing each subpath by a 2-vertex. The construction of G'** ***doesn't replace*** **our original construction of G. We still require G to obtain the proper edge emanating from or terminating in a 2-vertex. After obtaining a hamilton circuit in $G_{m''}'$ ($D_{m''}'$), we replace each 2-vertex by the subpath it represents to obtain a hamilton circuit in $G_{m''}$ ($D_{m''}$). In general, $G_{m''}'$ has minimum degree 3, while $\delta^+(D_{m''}') \geq 2$, $\delta^-(D_{m''}') \geq 2$. Henceforth, for simplicity, let n be the number of vertices in both the original graph and its contracted graph. For all graphs and directed graphs other than $D_{m''}'$, we construct $h_0$ by randomly choosing vertices from the balanced, binary search tree obtained from S = {1, 2, ... , n}, deleting entries from S after they are chosen after which we rebalance the search tree. Let 1 be the initial vertex of ORD($h_0$) and define ORD(1) = 1. If $v_2$ is the second vertex chosen, let ORD($v_2$) = 2, etc. ... . As we construct $h_0$, we place each successive number, $v_i$, along with its ordinal number, ORD($v_i$) ,in a balanced binary search tree, $H_0$, in** ***which the search key is ORD($v_i$).*** **Thus, if ORD(a) < ORD(b), then a occurs before b, going in a clockwise manner, starting at 1, around $h_0$. We, simultaneously, construct a balanced, binary search tree, ORD($H_0$) , in which the numbers from 1 through n occur in sequential order, in which each integer is followed by its ordinal number with respect to $h_0$.** ***Here the search key is the numerical value of each number from 1 through n.*** **This allows us to access the ordinal value of any given vertex on $h_0$ in at most O(logn) running time. As we construct $h_0$, we check each new arc in the pseudo-hamilton circuit, $H_0$, to see if it is a pseudo-arc or an arc of G. If it is a pseudo-arc vertex, v, we place the edge, [v, v"] , which contains the arc (v, v"), in a balanced, binary search tree called PSEUDO which has the following properties: One branch of PSEUDO is arranged so that if v < v', then [v, v"] precedes [v', v"'] on the tree. The second part of the tree follows the same rule for the second elements of the ordered pairs. Iterations of our algorithm replace $H_0$ by successive pseudo-hamilton circuits, $H_i$ (i = 1,2,3, ... ) where $H_i$ replaces $H_{i-1}$. Let j be the number of successes and i the number of failures in ITER = i + j iterations of the algorithm. From Lemmas 2 and 3 and the Law of Large Numbers, for very large values of n, the number of times we succeed minus the number of times we fail is at least $(\frac{12}{16} - \frac{3}{16})j = (\frac{9}{16})j$. In** general, if we have two or more $H_i$-admissible permutations, we do not use one which is the inverse of $s_{i-1}$ where $H_{i-1}\sigma_{i-1} = H_i$. For a large number of iterations, the net number of *successful* iterations in the first phase of the algorithm is at least

$$(.5625)(2n\log n) = 1.125n\log n$$

iterations. This number is large enough so that we can almost always go through each vertex of G' at least once. This implies that we almost always will obtain a hamilton path. In the second phase of the algorithm, we apply another 2nlogn iterations to almost always obtain a hamilton circuit in G'. For i ≥ 0, the iterations of our algorithm replace pseudo-arcs on $H_i$ by edges in G'. Let a be a pseudo-arc vertex of $H_0$. Then there are always *at least* six possibilities for choices of edges with respect to a as shown in Theorem E:

(1) [a, $H_0$(b)] and [b, $H_0$(c)] belong to G' - $H_0$;

(2) [a, $H_0$(b')] and [b', $H_0$(c')] belong to G' - $H_0$;

(3) [a, $H_0$(b)] and [$H_0$(a), d] belong to G' - $H_0$;

(4) [a, $H_0$(b')] and [$H_0$(a), d] belong to G' - $H_0$;

(5) [a, $H_0$(b)] and [$H_0$(a), d'] belong to G' - $H_0$;

(6) [a, $H_0$(b')] and [$H_0$(a), d'] belong to G' - $H_0$.

We now consider the probability of success during an iteration in $G_{m''}$' assuming the hypotheses of Theorem ABKS. Consider the worst possible cases. In what follows, we assume that $H_0$(a) and b are both arc vertices. Furthermore, assume that each of b and b' has two edges incident to it lying on $H_0$, while, since a is a pseudo-arc vertex, each of a and $H_0$(a) has precisely one edge incident to it lying on $H_0$. We now randomly choose up to logn edges incident to a and up to logn incident to $H_0$(a). (In any case, since every vertex of $G_{m''}$' is at least of degree three, we can always test at the minimum the six cases given above.) We systematically test pairs of edges – one incident to a, the other incident to b, b', or $H_0$(a) – for $H_0$-admissibility as pseudo-3-cycles. If we obtain an $h_0$-admissible pseudo-3-cycle, say $s_0$, we obtain $h_1 = h_0 s_0$ which implies that $H_1 = H_0\sigma_0$. Suppose we cannot obtain an $h_0$-admissible pseudo-3-cycle and we have two or more pseudo-arc-vertices on $H_0$, say a and v. Then we randomly choose up to logn edges incident to a and another set of up to logn edges incident to v. Using pairs of edges – one edge incident to a, the other, incident to v – we test for $h_0$-admissible pairs of POTDT. (We again note that $h_0$-admissibility occurs if and only if the two edges intersect.) If s is an $h_0$-admissible POTDT, then $h_1 = h_0 s_0$ implying that $H_1 = H_0\sigma_0$. Before going on, let DIFF = $e_s$ - $a_s$ where $e_s$ is the number of edges in G' - $H_0$ associated with an $h_0$-admissible permutation, s, and $a_s$ is the number of arc vertices in s. If, after a full search of possibilities for $h_0$-admissibility, we obtain two or more successful outcomes, we first check to see which among the cases yields the largest value of DIFF. If more than one such case exits, we pick one of them, say [a, $H_0$($b_j$)] (j ≤ [logn]), which is contained in the largest number of admissible permutations obtained. Whether we succeed or fail to obtain an $h_0$-admissible permutation, at the end of the iteration, we do a rotation of the following kind: Given a pseudo-arc vertex of $H_0$, say v, we randomly choose an edge incident to v which doesn't lie on $H_0$, say [v, v']. Assume that the following subpath lies on $H_0$:

$$S = [v, v_1, v_2, \dots, v_r, v']$$

Using a rotation with respect to $v_1$ and v', S is transformed into

$$S' = [v, v', v_r, \dots, v_1]$$

Note that the pseudo-arc vertex, v, is now an arc vertex, while the (possibly arc) vertex, $v_1$, in most cases becomes a pseudo-arc vertex. (The procedure of randomly choosing a rotation after an iteration – successful or not – is especially useful if we wish to obtain a hamilton circuit in a non-random graph. It insures that we can randomly go through all of the vertices of the graph. Otherwise, it is possible that we may go through a loop in which we only go through a fixed subset of V.) The use of abbreviations shortens the running time of the algorithm. We test (a b c) for $H_0$-admissibility in the following way: Using ORD($H_0$), we find the ordinal values of a, b and c with respect to $H_0$. Then (a b c) is $H_0$-admissible if and only if ORD(a), ORD(b) and ORD(c) occur in a clockwise manner. If (a b c) is $H_0$-admissible, we form a balanced, binary search tree called $A_i$ where

$$A_1 = [ORD(a), H_0(ORD(b)), \dots, ORD(c) H_0(ORD(a)), \dots, ORD(b) H_0(ORD(c))]$$

If (a b c) is $H_0$-admissible, we form a queue called BACKTRACK in which we place (a c b). Since

$$H_0 (a\ b\ c) = H_1,\ H_0 = H_1 (a\ b\ c)^{-1} = H_1 (a\ c\ b).$$

***Note. When we backtrack, it is possible that more than one of the edges in $S = \{[a, H_1(c)], [b, H_1(a)], [c, H_1(b)]\}$ is a pseudo-arc.***

As we proceed in the construction of the abbreviations $A_i$ ($i = 1, 2, \dots, [n^{.5}]$), we use the following rules:

Let v be a random element of V.

(1) If the successor to ORD(v) doesn't explicitly occur in $A_i$, then $H_i$(ORD(v)) = $H_0$(ORD(v));

(2) otherwise, $H_i$(ORD(v)) = the successor of ORD(v) on $A_i$.

To clarify the construction of $A_1$ and $A_2$, we use the following example:

Let

$$H_0 = (1\ 14\ 8\ 4\ 3\ 12\ 7\ 13\ 10\ 6\ 11\ 5\ 15\ 9\ 2)$$
$$= (ORD(1)ORD(2) \dots ORD(11)ORD(12)ORD(13)ORD(14)ORD(15))$$

From the latter, we note that the orders of the elements of the permutation are a subset of the natural numbers. It follows that in a first rotation, those ordinal numbers which do not appear in $A_1$ and which are assumed to have negative signs in front of them must be a subset of {n, n-1, n-2, ... , i, i-1, ...}. Thus, we come to the following rule for rotations:

(1) During a rotation with respect to v and $v_i$, [v, $v_i$] becomes an arc of $H_0$, while [$v_{i-1}$, $v_{i+1}$] may or may not change its orientation. All other arcs or pseudo-arcs formed have the same designation (pseudo-arc or arc) as they had previously. More simply, if $(v_c, v_{c+1})$ was a pseudo-arc before the rotation, then $(v_{c+1}, v_c)$ would be a pseudo-arc after the rotation. Similarly, if $(v_c, v_{c+1})$ was an arc before the rotation, then $(v_{c+1}, v_c)$ is an arc after the rotation. The reasoning is straight-forward here: Both arcs come from the same edge which either lies on $H_0$ or doesn't lie there.

(2) If an edge of $H_0$, say $[v_c, v_{c+1}]$, lies in PSEUDO and has elements which do not explicitly occur in an abbreviation, then if a "-" precedes each of element, $v_{c+1}$ precedes $v_c$ and is a pseudo-arc vertex. If no sign precedes each of them, then $v_c$ is a pseudo-arc vertex. Suppose v lies in an abbreviation and we can infer that v" is its successor on $H_0$. Then if [v, v"] or [v.'s, v] lies in PSEUDO, from (1) we know that v is a pseudo-arc vertex. The same is true if both vertices lie in an abbreviation.

Continuing with our example,

ORD(1) = 1, ORD(2) = 14, ORD(3) = 8, ORD (4) = 4,
ORD(5) = 3, ORD(6) = 12, ORD(7) = 7, ORD(8) = 13,
ORD(9) = 10, ORD(10) = 6, ORD(11) = 11, ORD(12) = 5,
ORD(13) = 15, ORD(14) = 9, ORD(15) = 2

Suppose $s_1$ = (1 4 7). Since ORD(1) = 1, ORD(4) = 4, ORD(7) = 7,
If 1, 4, 7 traverse $H_0$ in a clockwise manner, $s_1$ is $H_0$-admissible. Using

$$H_0(\text{ORD}(1)) = \text{ORD}(2),$$
$$H_0(\text{ORD}(4)) = \text{ORD}(5),$$
$$H_0(\text{ORD}(7)) = \text{ORD}(8),$$

$$A_1 = \{\text{ORD}(a)H_0(\text{ORD}(b))\ldots\text{ORD}(c)H_0(\text{ORD}(a))\ldots\text{ORD}(b)H_0(\text{ORD}(c)\ldots\}$$
$$= \{\text{ORD}(1)\text{ORD}(5)\ldots\text{ORD}(7)\text{ORD}(2)\ldots\text{ORD}(4)\text{ORD}(8)\ldots\}$$

Henceforth, we simplify notation by using ordinal numbers in abbreviations. Assume now that 7 is a pseudo-arc vertex and that we ant to construct a rotation out of 7. W.l.o.g., let [7, 10] belong to $G_{m*}'$. ORD(7) = 7, ORD(10) = 9. Then the rotation with respect to ORD(7) and ORD(9) transforms $A_1$ into $A_2$ in the following way:

$A_2$ = {1 5 ... <u>7 9</u> 8 4 --- <u>2 10</u> ...}

The dashes between 4 and 2 indicate that the ordinal numbers are consecutively decreasing in value. Thus, [7 2 3 4 8 9 10] in $A_1$ becomes [7 9 8 4 3 2 10]. The important thing to note is rotations never increase the number of pseudo-arc vertices in $H_1$. First, 7 is no longer a pseudo-arc vertex since it is followed by 9 where [ORD(7), ORD(9)] = [7, 10] is an edge of $G_{m''}'$. On the other hand, [ORD(2), ORD(10)] is generally not an edge of $G_{m''}'$, although if rotations are done often enough it may in a particular case belong to $G_{m''}'$. Now suppose [ORD(4), ORD(8)] belongs to $G_{m''}'$. Then

$$[\text{ORD}(8), \text{ORD}(4)] = [\text{ORD}(4), \text{ORD}(8)]$$

also belongs to $G_{m''}'$, while if ORD(3) is a pseudo-arc vertex, then [ORD(3), ORD(4)] = [ORD(4), ORD(3)] doesn't belong to the graph. It follows that if ORD(3) is a pseudo-arc vertex of $A_1$, then ORD(4) is a pseudo-arc vertex of $A_2$. Now we want to construct $A_3$. Let us assume that ORD(2) = 14 is a pseudo-arc vertex. W.l.o.g., let $s_2$ = (14 9 12).

14 = ORD(2), 9 = ORD(14), 12 = ORD(6)

We place the respective ordinal numbers in $A_2$ underlined and in italics to see if they occur in a clockwise manner (going from left to right and starting at 1 again if necessary).

$A_2$ = (1 5 *6* ... 7 9 8 4 --- *2* 10 ... *14* ...)

Going from left to right, we obtain 6, 2, 14. The numbers occur in a clockwise manner in the circle defined by $H_2$. Thus, $s_2$ is $H_2$-admissible. Before we can construct $A_3$ to represent $H_3$, we must obtain $H_2$(ORD(2), $H_2$(ORD(14)), $H_2$(ORD(6)). The successor of ORD(2) in $A_2$ is ORD(10). On the other hand, the successor of ORD(14) ***doesn't explicitly occur*** in $A_2$. Therefore, it is its successor in $H_0$, namely, ORD(15). Finally, consider ORD(6). The successor of ORD(6) doesn't explicitly occur in $A_2$. Therefore, its successor is its successor in $H_0$, namely, ORD(7). Thus,

$6 \rightarrow 2 \rightarrow 10,\ 2 \rightarrow 14 \rightarrow 15,\ 14 \rightarrow 6 \rightarrow 7$

yielding

$A_3$ = (1 5 *6 10* ... *14 7* 9 8 4 --- *2 15*)

Now let 12 = ORD(6) and 3 = ORD(5) be pseudo-arc vertices. Let

$s_3$ = (ORD(5) ORD(13))(ORD(6) ORD(2)) = (3 5)(12 14)

be a product of two disjoint pseudo-2-cycles (POTDTC) which we wish to test for $H_3$-admissibility. From

$A_3$ = (1 *5* *6* 10 ... *13* 14 7 9 8 4 --- *2* 15)

we see that [5 13] intersects [6 2] in the circle $H_3$. Therefore, $s_3$ is $H_3$-admissible. In this case, $A_4$ is of the form

(a$H_3$(b) ... d$H_3$(c) ... b$H_3$(a) ... c$H_3$(d) ...)

Continuing,

$H_3$(ORD(5)) = ORD(6), $H_3$(ORD(6)) = ORD(10),
$H_3$(ORD(13)) = ORD(14), $H_3$(ORD(2)) = ORD(15)

We thus obtain

$A_4$ = (1 5 14 7 9 8 4 --- 2 10 ... 13 6 15)

Before we apply a new permutation in an iteration, we check BACKTRACK to see if the new permutation, $s_i$, is at the end of the queue. If it is, we skip it and go on with our search. If $G'$ is a directed graph, and the only $H_i$-admissible permutation is at the end of the queue on BACKTRACK, we use it to backtrack to $H_{i-1}$. On the other hand, if $G'$ is a graph and the only $H_i$-admissible permutation is at the end of the queue in BACKTRACK, instead of backtracking, we construct a rotation to obtain $H_{i+1}$ and continue with the algorithm. After we apply a new permutation, (a b c) or (a c)(b d) to an abbreviation, we place its inverse ((a c b) or (a c)(b d)) in BACKTRACK. In general, whether we succeed or fail in an iteration, if [v, v""] changes from a pseudo-arc on $H_i$ to an arc, [v, $H_i$ (v")], on $H_{i+1}$, we delete [v, v"] from PSEUDO. On the other hand, if [v', v""] is a new pseudo-arc on $H_{i+1}$, we place [v', v""] in PSEUDO. Going back to our example, if we fail in an iteration applied to $A_4$ and ORD(13) is a pseudo-arc vertex, we construct a new rotation using an edge of $G_{m''}$' not on $H_4$, say [ORD(13), ORD(7)]. After the rotation, [ORD(6), ORD(9)] will most often become a pseudo-arc, while [ORD(13), ORD(7)] will become an arc of $H_5$. If [ORD(6), ORD(9)] is an arc of $H_5$, we obtain a pseudo-arc from PSEUDO. If PSEUDO contains no pseudo-arc, then $H_5$ is a hamilton circuit. In general, if G is not a directed graph and we don't have to backtrack, the number of arcs in PSEUDO is a monotonically decreasing function which approaches 1 as we go through all of the vertices in G. If G is a directed graph, we generally must backtrack; in some cases when we backtrack, we increase the number of pseudo-arcs by 1 or 2. When we reach the abbreviation $A_{[n^5]}$, we construct $H_{[n^5]}$ and ORD($H_{[n^5]}$) using $A_{[n^5]}$ and $H_0$. We then delete $H_0$ and $A_{[n^5]}$; we next use $H_{[n^5]}$ to construct abbreviations $A_{[n^5]+1}$, $A_{[n^5]+2}$, ... , $A_{2[n^5]}$. Using $H_{[n^5]}$ and $A_{2[n^5]}$, we construct $H_{2[n^5]}$, and then delete $H_{[n^5]}$ and $A_{2[n^5]}$. This procedure continues throughout the algorithm. From Theorems A, B and E, the probability of an iteration yielding at least two $H_i$-admissible permutations is at least $\frac{3}{4}$. For very large n, if we use 2nlogn iterations - where each iteration uses up to 2(logn)$^2$ edges randomly chosen to obtain an admissible permutation - we almost always successfully pass through every vertex in V. If there are fewer than logn edges incident to a vertex, a, we may use each edge in G' - $H_i$ incident to a. If there are two pseudo-arc vertices, a and b, with which we try to construct a pseudo-POTDT, we may use up to logn edges incident to each of them in the constructions. Using 2nlogn iterations, the probability that we will be able to obtain a hamilton path approaches 1 as n to inf. After obtaining a hamilton path, say H'P, which contains only one pseudo-arc vertex, a hamilton circuit, H'C, is obtained using another 2nlogn iterations: Since G' is a random graph, at some point we obtain an $H_i$-admissible 3-cycle, say (p q r), such that each edge in

$$S = \{[p, H_i(q)], [q, H_i(r)], [r, H_i(p)]\}$$

lies in G' - $H_i$. In order to complete the algorithm, we must replace each 2-vertex by the subpath which it represents to obtain a hamilton circuit in $G_{m*}$. Let $A_{h'p+j}$ be the abbreviation in which we obtain H'C, while $H_{h'p}$ is hamilton path associated with $A_{h'p+j}$. Using $H_{h'p}$ and the sign in front of a 2-vertex, we can determine the correct orientation of the subpath in the hamilton circuit, HC, of $G_{m''}$. Alternately, given the 2-vertex, vv', we can determine the correct orientation by determining which vertex the edge *entering* vv' is incident to. By construction, the edge *emanating* from vv' is incident to the other vertex. We call the algorithm just described, ALGORITHM G.
The problem in applying the algorithm to a directed graph is that we can't use rotations in it. Thus, we actually have to backtrack. If backtracking we always advance the index of $H_i$.

**Thus, $H_i$ (a c b) = $H_{i+1}$ = $H_{i-1}$. Since (a c)(b d) is its own inverse, if necessary, we place (a c)(b d) in BACKTRACK. We always *assume* that if (a b c) is $H_i$-admissible, then (a c b) is $H_{i+1}$-admissible. However, if $H_{i+1}$ has *fewer* pseudo-arc vertices than $H_i$,**

$$S_{i+1} = \{[a, H_{i+1}(c)], [c, H_{i+1}(b)], [b, H_{i+!}(a)]\}$$

**contains *fewer than two arcs*. Although we use the inverse of $s_{i+1}$ ((a c b) or (a c)(b d)) in all cases, we should keep this fact in mind. On the other hand, if $H_i$ and $H_{i+1}$ have the same number of pseudo-arc vertices, then $S_{i+1}$ always contains at least two arcs. This is why the probability must be greater than $\frac{1}{2}$ that we have at least *two* admissible permutations in an iteration: if we have only *one*, say (a c b), it may well be that we have to use it to backtrack when we're working with a directed graph. In the case of a directed graph, $D_{m''}$, its contracted graph, $D_{m''}'$, has at least two arcs entering and two leaving each vertex in V. We have no trouble defining the orientation of any subpath $S_\beta$. It has fixed orientation throughout the algorithm. Thus, if a 2-vertex is vv', its initial vertex is v and its terminal vertex is v'. Also, since we don't use rotations, PSEUDO consists of pseudo-arc vertices – not pseudo-arcs. Otherwise, the algorithm is the same as ALGORITHM G. WE call the algorithm for directed graphs ALGORITHM D. From Theorem E, the probability that G has at least two $H_i$-admissible permutations approaches $\frac{13}{16}$ as n to inf. Correspondingly, the probability for failure approaches $\frac{3}{16}$. It follows that the net number of successful iterations is at least $\frac{1}{2}$ ITER where ITER is the number of iterations. Assuming that n is very large, we again use 2nlogn iterations in order to successfully go through each vertex in $D_{m''}'$. We then require another 2nlogn iterations to obtain a hamilton circuit.**

## IV. PROBABILITY OF SUCCESS.

**In any directed graph, D, considered here,**

$$\delta^+(D' - H_i) \geq 1, \ \delta^-(D' - H_i) \geq 1$$

**i = 0,1,2, ... . Furthermore, if a is a pseudo-arc vertex,**
**$\delta^+(a) \geq 2$, $\delta^-(a) \geq 1$ in D' - $H_i$, while $\delta^+(H_i(a)) \geq 1$, $\delta^-(H_i(a)) \geq 2$**
**(i = 0,1,2, ... ).**
**Alternately, if G is a graph, then $\delta(G' - H_i) \geq 1 (i = 0,1, ...)$**
**Furthermore, if a is a pseudo-vertex of $H_i$, then**

$$\delta(a) \geq 2 \text{ in } G' - H_i, \ \delta(H_i(a)) \geq 2 \text{ in } G' - H_i$$

**Therefore, the probability of success when searching in depth in up to $(\log n)^2$ trials in each iteration is at least $\frac{3}{4} = \frac{12}{16}$. By the Law of Large Numbers, it follows that the number of successes in 2nlogn iterations approaches at least 1.5nlogn, while the number of failures is at**

most .5logn. It follows that the net number of successes approaches at least nlogn. Thus, if v is an arbitrary vertex in G', the probability that a successful iteration passes through v is at least $\frac{1}{n-1}$. Thus, the probability that a successful iteration passes through *each* vertex of G' is at least

$$1 - (1 - \frac{1}{n-1})^{n\log n} \rightarrow 1 - e^{-\log n} \rightarrow 1 - \frac{1}{n} \rightarrow 1$$

as n to inf. It follows that the probability of obtaining a hamilton path in G' (and therefore in G) approaches 1 as n to inf. Assume now that $H_P$ is a hamilton path in G'. Let i ≤ 2nlogn. Then the probability that an $H_{P+i}$-admissible 3-cycle has all of its edges in G' - $H_{P+i}$ is $\frac{1}{n-1}$. It follows that the probability of obtaining a hamilton circuit in G' (and thus in G) is

$$1 - (1 - \frac{1}{n-1})^{n\log n} \rightarrow 1 - e^{-\log n} \rightarrow 1 - \frac{1}{n} \rightarrow 1$$

as n to inf, concluding the proof.

## V. RUNNING TIME.

The steps that follow are required for running through Algorithm G and Algorithm D. For graphs and directed graphs which have a minimal degree greater than 2, steps 2* and 3* my be omitted.

(1) Constructing G.

(2*) Searching for vertices of degree 2 in $G_{m^*}$ and $D_{m^*}$ and constructing corresponding subpaths, say $S_\gamma$.

(3*) Replacing each $S_\gamma$ by a 2-vertex and constructing the respective graph and directed graph, $G_{m''}'$, $D_{m''}'$.

(4) Constructing $h_0$, $H_0$ and ORD($H_0$), $H_0'$, ORD($H_0'$).

(5) Constructing and working with PSEUDO.

(6) Systematically testing up to logn edges out of a pseudo-vertex, a. Then given a chosen edge, [a, $H_i$(b)],

(a) in up to [logn] cases, checking whether edges [b, $H_i$($c_j$)] intersect [a, $H_i$(b)];

(b) in up to [logn] cases, checking whether edges [$H_i$(a), $d_k$] intersect [a, $H_i$(b)].

(6*) Systematically testing up to [logn] edges respectively out of a pseudo-arc vertex, a*, and a pseudo-arc vertex, b*, to obtain all possible pairs which intersect.

(7) Constructing $A_i$'s and $H_j$'s for

$$j = [in^{.5}] \ (i = 0, 1, \ldots, [n\log n])$$

(8) Constructing a rotation out of a pseudo-arc vertex of a graph.

(9) Constructing BACKTRACK.

We can eliminate (9) for Algorithm D. In addition, if an abbreviation, $A_i$, has j entries, it requires O(j) running time to apply a rotation to it. Thus, the running time of Algorithm G is not greater than that of Algorithm D.

(10) Replacing 2-vertices in $H_C'$ by the subpaths they represent and constructing $H_C$.

(1) Constructing G.

We first construct a balanced, binary search tree, $G_{TREE}$, containing 1, 2, ... , n as distinct branches. then we randomly choose arcs from V X V, placing the terminal point on each arc on the correct branch. It takes at most $(\log n)^2$ time to randomly pick an arc from V X V. It takes a further logn time to place each arc chosen on G. Since the number of edges chosen is never greater than O(nlogn), the running time necessary to construct G is not greater than $O(n\log^3 n)$.

(2*) In Algorithm G, it requires O(nlogn) time to discover each vertex of degree 2. From the Corollary to Theorem B, as n to inf, the number of vertices of degree 2 is at most O(logn). Given a fixed vertex of degree 2, the probability that an adjacent vertex is of degree 2 is at most $\frac{2\log n}{n}$. There are at most $\binom{\log n}{2}$ pairs of vertices of degree 2. Thus, the probability that *no* pair of vertices of degree 2 are adjacent is at least

$$(1 - \frac{2\log n}{n})^{(\log n)^2} \approx \exp(-\frac{2(\log n)^3}{n}) \;\to 1$$

as n to inf. It follows that if n is very large we almost always have at most one vertex of degree 2 in a subpath. Therefore, it takes at most $O(n(\log n)^2)$ running time to construct all subpaths, $S_\alpha$, in G and place them in a balanced binary search tree, TWOPATH. From the Corollary to Theorem C, the probability that a unique arc emanates from, or terminates more than $O((\log n)^2)$ vertices approaches 0 as n to inf. Thus, in Algorithm D, it requires at most $O(n(\log n)^3)$ running time to obtain all subpaths, $S_\beta$, and place them in a balanced, binary search tree, DTWOPATH.

(3*) Given Algorithm G or Algorithm D, it takes at most $O(n(\log n)^3)$ running time to construct a duplicate graph of G in which each 2-subpath of TWOPATH (DTWOPATH) is replaced by a 2-vertex. This involves deleting the interior edges and vertices of each subpath which takes $O((\log n)^3)$ r. t.. In particular, if $v_1v_2$ is a 2-vertex, the branch of $v_1$, $B_{v_1}$, and the branch of $v_2$, $B_{v_2}$, each becomes a distinct branch of $v_1v_2$. The latter requires adding each vertex of $B_i$ (i = 1,2) to a new branch, $B_{i'}$, of $v_1v_2$ and then deleting the entries of $B_i$ as well as $v_i$. We name the new graph (directed graph) $G_{m''}'$ ($D_{m''}'$). For large n, the number of 2-paths and 2-vertices is no greater than $O((\log n)^2)$. It requires at most $O(n(\log n)^3)$ r. t. to construct a duplicate of G. It thus requires at most $O((\log n)^3)$ r. t. to make deletions and insertions in the construction of $G_{m''}'$ ($D_{m''}'$). It thus requires at most $O(n(\log n)^3)$ r. t. to construct $G_{m''}'$ ($D_{m''}'$).

(4) Constructing $h_0$, $H_0$ and ORD($H_0$), $H_0'$ and ORD($H_0'$).

In cases other than $D_{m''}$, given that V is on a balanced, binary search tree, it requires logn n time for each of the following operations: pick a random number from the tree, delete that number, rebalance the tree. If when constructing $H_0$, we assume that ORD($S_\alpha$) or ORD($S_\beta$) has a single ordinal value, then we can simultaneous construct $H_0$ and $H_0'$. Furthermore, since all we requires is ORD($H_0'$), we can construct it directly.(We note that if $\delta(G) > 2$, then
$H_0 = H_0'$.) Thus, it requires O(nlogn) time to construct each of $h_0$, $H_0$ and ORD($H_0$).

(5) Constructing and working with PSEUDO.

As we construct $h_0$, we add at most n vertices in a directed graph and 2n vertices from n pseudo-arcs in a graph. It requires at most 2logn time to do each of the following: add a vertex of pseudo-arc, delete a vertex or pseudo-arc, rebalance the binary search tree. Accordingly, to do this in 4nlogn iterations, it requires at most O(nlogn) running time. For an arbitrary pseudo-3-cycle, we must add and/or delete at most 12 vertices. For an arbitrary POTDT, we must add and/or delete at most 16 vertices. This requires no more than O(logn) running time. The number of iterations in which we do this is [4nlogn]. Thus the total running time for working with PSEUDO is at most $O(n(\log n)^2)$.

(6)-(6*)-(7).

There is a counter indicating the number of edges incident to each vertex in a balanced, binary search tree represent a graph. Similarly, in a directed graph, there is a counter indicating the number of arcs emanating from, or terminating in each vertex. Thus, the running time for obtaining at most $2([\log n])^2$ edges is at most $2([\log n])^3$. Therefore, the running time for choosing arcs in [4nlogn] iterations is at most $O(n(\log n)^4)$. If (a b c) is $H_0'$-admissible, we place the ordered set

$$\{ORD(a),ORD(H_0'(b)),..,ORD(c),ORD(H_0'(a)),..,ORD(b).ORD(H_0'(c)),..\}$$

in the balanced, binary search tree, $A_1$. Since ORD($H_0'$) contains the ordinal value of each vertex of V, it requires at most 6logn running time to locate the ordinal values with respect to $H_0'$ of each of these six vertices and place them in $A_1$. (With an $H_0'$-admissible POTDT, it would require at most 8logn r. t.) Thus, testing up to $2([\log n])^2$ edges requires at most $16([\log n])^3$ r. t. Using the ordinal values of another set of at most eight vertices and placing them properly among the elements of $A_1$ to form $A_2$ requires $2(16[\log n])^3$ r. t. . It follows that to construct each $A_i$ ($i = 1,2, \ldots ,[n^{.5}]$) requires

$$16([\log n])^2 \sum_{i=1}^{i=[n^{.5}]} i < 8n([\log n])^2$$

r. t., implying that these operations applied $[4n^{.5}\log n]$ times requires require $O(n^{1.5}(\log n)^4)$ r. t. We can do each of the following in O(nlogn) r. t.:

(i) reconstruct $H_{[n^{.5}]}'$ from $A_{[n^{.5}]}$ and ORD($H_0'$);

(ii) delete $H_0'$ and construct ORD($H_{[n^{.5}]}'$)

**Since we require 4[n$^{.5}$ logn] iterations to complete the algorithm, and we go through [n$^{.5}$] iterations before replacing $H_{(i-1)[n^{.5}]}$' by**

**$H_{i[n^{.5}]}$', it requires O(n$^{1.5}$ (logn)$^4$) r. t. to complete this portion of the algorithm.**

**(8) Constructing a rotation out of a pseudo-arc.**

**We can do the rotation on an abbreviation containing i vertices. It requires no more than O(logi) time to locate the two vertices defining the rotation. We then require no more than O(logi) time to make changes on the abbreviation. We make a rotation at the end of each iteration. Thus, using the computations given in the (6)-(6*)-(7), it requires at most O(n$^{1.5}$ (logn)$^4$) to make of all the rotations on a graph using Algorithm G.**

**(9) Constructing BACKTRACK.**

**Assume that $H_{i+1}$' = $H_i$'(a b c) or $H_i$'(a c)(b d) is a pseudo-hamilton circuit in a directed graph D'. Suppose that we fail to obtain an $H_{i+1}$'-admissible permutation during an iteration. Then we must obtain either**

**$S_{iii}$ = {a, c, b} or $S_{iv}$ = {a, c, b, d}**

**from the near end of the queue BACKTRACK to apply it $A_{i+1}$ to obtain $A_{i+2}$ = $A_i$. We then delete the elements in $S_{iii}$ or $S_{iv}$ from BACKTRACK and subtract 1 from the counter. Generalizing, suppose that S is the set of points of an $H_i$'-admissible permutation and that $A_i$ contains at most 4i points. Even if this deletion occurs for ever iteration, it requires at most 16logn time to obtain each ordinal value from ORD($H_i$') and 16i time to find them and delete them from BACKTRACK. It thus requires at most O(n$^{1.5}$ (logn)$^2$) to On the other hand, it requires at most 16nlogn time to place permutations in BACKTRACK. It therefore requires at most O(n$^{1.5}$ (logn)$^2$) to work on BACKTRACK.**

**(10) Replacing 2-vertices in $H_C$' by subpaths to form $H_C$.**

**There are at most O([logn])$^2$) subpaths. Furthermore, there are at most four vertices on a subpath. It thus requires at most 4logn r. t. to locate each of them and another 4(logn)$^2$ r. t. to replace a 2-vertex on $H_C$' by the corresponding subpath to form $H_C$. The total operation therefore requires at most O((logn)$^4$) rt. .**

**Reviewing (1)-(10), the running time of either algorithm is at most O(n$^{1.5}$ (logn)$^4$) r. t. .**

## V. FURTHER RESULTS.

**Henceforth, given a graph or directed graph, G, the edges of the pseudo-hamilton circuit, $H_0$', lie in $K_n$ - G'. Also, any of the following will be denoted by as an $H_i$-admissible permutation: $H_i$-admissible 3-cycle, $H_i$-admissible pseudo-3-cycle, $H_i$-admissible POTD 2-cycles, $H_i$-admissible POTD pseudo-2-cycles. For simplicity, an $H_i$-admissible POTD pseudo-2-cycles may have one or two pseudo-2-cycles.**

**THEOREM F. If G is a random directed graph, as n to inf, there always exists a hamilton circuit, $H_C$ ', in G' obtainable from an arbitrary pseudo-hamilton circuit, $H_0$ ', by sequentially using only $H_i$ '-admissible permutations without backtracking.**

**PROOF. Algorithm G proves the theorem except for backtracking. But, as shown in section I,**

$$h_C ' = (h_0 ') \prod_{i=1}^{i=c+p} \sigma_i$$

**where p stands for the number of permutations which were used in backtracking. But each of these permutations changed $H_i$ ' to $H_{i-1}$ '. Thus, each of these permutations was the inverse of the preceding one. It follows that, starting at $H_C$ ' and removing all of the permutations which were erased by backtracking, we are left by a product of 3-cycles each of which is $H_{\alpha_i}$ '-admissible (i=0,1,2, ... ,N) where N is less than 4nlogn. This product yields $H_C$ '.**

**Using the result of Theorem F, we can use the theorem in quantum computing by eliminating branches which require backtracking. Theoretically, we could obtain all hamilton circuits of G. Whether such a computer can ever be built is another question.**

**Given the statement preceding Theorem F, we may assume that $H_0$ ' contains only pseudo-arc vertices. Before proving Theorem G, we prove the following lemmas:**

**LEMMA 1q. Assume that we obtain $H_i$ by applying an $H_{i-1}$-admissible permutation to $H_{i-1}$ (i = 1,2, ... ). Let $\sigma_i = (H_i)^{-1} H_C$ where $H_C$ is a hamilton circuit of a graph, G. Then the pseudo-arc vertices (the initial vertices of the pseudo-arcs) in PSEUDO are the same as the points moved by $\sigma_i$.**
**PROOF. All of the edges of $H_0$ lie in $K_n$ - G. Thus, the only edges lying on $H_i$ (i = 1,2, ... ) are edges of $H_C$. It follows that for any value of i, if $\sigma_i = (H_i)^{-1} H_C$, any identity element of the permutation $\sigma_i$ corresponds to an edge, say e, of $H_C$ which lies on $H_i$. It follows that e is an edge of $H_i$ whose initial vertex is an arc vertex. On the other hand, all edges of $H_C$ which have initial vertices moved by $\sigma_i$ are pseudo-arc vertices of $H_i$. It follows that |PSEUDO| = $|\sigma_i|$.**

**LEMMA 2q. Let G be a graph containing a hamilton circuit $H_C$. Let $H_i$ (i = 0,1,2, ... ,) be a pseudo-hamilton circuit of G obtained from successively applying admissible permutations to $H_0$, $H_1$, ... , $H_{i-1}$. Assume $\sigma_i = (H_i)^{-1} H_C$. Suppose all of the points of a disjoint cycle of $\sigma_i$, say C, do not traverse $H_i$ in a counter-clockwise manner, | C | > 2, and either**

**(i) going in a clockwise manner, three consecutive points of C, say a, b, c traverse $H_i$ in a clockwise manner,**

**(ii) three consecutive points of C, say c, a, b traverse $H_i$ in a clockwise manner, while c, b, a traverse $H_i'$ in a counter-clockwise manner.**
**Then (a b c) is an $H_i$-admissible permutation.**

**PROOF.**
**If a, b, c traverse $H_i$ in a clockwise manner, then**

**(i) [a, $H_i$(b)] and [b, $H_i$(c)] in the circle containing the evenly- spaced vertices of $H_i$. It follows from Theorem A that (a b c) is an $H_i$-admissible permutation.**

**(ii) Let c, b, a traverse $H_i$ in a counterclockwise manner. It follows that a, b and c traverse $H_i$ in a clockwise manner with (a b) and (c a) arcs of C. Therefore, (a b c) is an $H_i'$-admissible permutation.**

**COROLLARY TO LEMMA 2q. Let $|s|$ denote the number of points moved by a permutation, s, of $S_n$. Then if $H_i = H_{i-1}$(a b c) and**
**$H_i \sigma_i = H_{i\ C}$, $|\sigma_i| \leq |\sigma_{i-1}|$ - 2.**

**PROOF. From LEMMA 1q, the vertices of $\sigma_i$ are precisely the elements of PSEUDO. Since (c a) and (a b) are both arcs of C, it follows that all of the vertices a, b, c are pseudo-arc vertices of $H_i$. Application of $H_{i-1}$ of (a b c) to form $H_{i-1}$ deletes at least two arcs from C as well as two pseudo-arcs or pseudo-arc vertices from $PSEUDO_{i-1}$. If one of the three arcs {(c a), (a b), (b c)} doesn't belong to $\sigma_{i-1}$, then $H_i \sigma_i = H_C$ where $|\sigma_i| = |\sigma_{i-1}|$ - 2. If all three arcs (a b), (b c) and (c a) belong to $\sigma_{i-1}$, then $|\sigma_i| = |\sigma_{i-1}|$ - 3.**

**LEMMA 3q. Let C be a disjoint cycle of $\sigma_i$ such that neither condition (i) nor condition (ii) of Lemma 2q holds for any three consecutive points. Then if (a c) is an arbitrary arc of C, there exists at least one arc, say (b d), of a cycle of $\sigma_i$ such that (a c)(b d) is an $H_i$-admissible permutation.**

**PROOF. If neither condition (i) nor condition (ii) holds for C, then we cannot obtain an $H_i$-admissible permutation of the form (a b c) using at least two adjacent arcs of C. (In what follows, we use the ordinal values in ORD($H_i$) of a, b and c. Since there is always a way of expressing the three ordinal values in the range [1, n], we assume that in each example this is the case. Given the previous statement, x > y means that ORD(x) > ORD(y).) If all sets of three consecutive points of C (say, for example, c, a, b) either traverse $H_i'$ in a counter-clockwise manner with**

**(i) c > a > b,**

**or else**

**(ii) c > a, a < b, b > c,**

**then (a b c) cannot be $H_i$-admissible. As an example, let**

**$H_i$' = (1 2 3 4 5 6 7 8 9 10), while c = 8, a = 5, b = 3. Then (a b c) = (5 3 8) which is not $H_i$-admissible. Since no set of three consecutive points of C form an $H_i$-admissible permutation, every arc, (a c), of C must interlace with an arc, (b d), belonging to some cycle (possibly C) of $\sigma_i$. For suppose that that is not the case. $H_i$ applied to the pseudo-2-cycle (or 2-cycle), (a c), yields the product of disjoint cycles**

$$(a\ H_i(c)\ \ldots)(c\ H_i(a)\ \ldots) = P_1 P_2$$

**where the points of $P_1$ and $P_2$ include all of the elements of {1, 2. ... ,n}. But $H_i \sigma_i = H_C$. Since there exists no arc of $\sigma_i$ interlacing with (a c), $P_1$ contains no point, p, which lies on an arc of $\sigma_i$ connecting p to some point outside of $P_1$. Thus, $P_1$ is a disjoint cycle of $H_C$ containing fewer than n points which is impossible since $H_C$, by definition, is a Hamilton circuit (an n-cycle).**

**EXAMPLE. Let $H_i$ = (1 2 3 4 5 6 7 8 9 10) with (a c) = (3 7). Then $H_i$(3 7) = (1 2 3 8 9 10)(4 5 6 7). Now let b = 4 and d = 9. Then $H_i$(3 7)(4 9) = H* = (1 2 _3_ 8 _9_ 5 6 _7_ _4_). Note the occurrence of 3, 7, 4, and 9 in H*. 3 – 9 – 7 – 4 : the points of (a c) interlace with those of (b d).**

**Note. A cycle which has an even number of points is an *odd permutation*, while one with an odd number of points is an *even permutation.* (odd)(odd) = even, while (odd)(even) = odd. Thus, a cycle containing n points when multiplied by (a c) can *never* yield another cycle containing n points.**

**LEMMA 4q. If $H_{i-1}$ is applied to an $H_{i-1}$-admissible permutation, (a c)(b d), to obtain an $H_i$-admissible permutation, then $|\sigma_i| \le |\sigma_{i-1}| - 2$.**

**PROOF. A product of two disjoint pseudo-2-cycles obtained from $\sigma_{i-1}$ contains four pseudo-arc vertices of $H_{i-1}$ and at least two arcs of $\sigma_{i-1}$. Thus, $|\sigma_i|$ is at least two less than $|\sigma_{i-1}|$.**

**COROLLARY TO LEMMA 4q. The number of pseudo-arcs or vertices in $PSEUDO_i$ is at least two fewer than the number in $PSEUDO_{i-1}$.**

**PROOF. From Lemma 1q, the number of pseudo-arcs or vertices in $PSEUDO_i$ equals the number of points moved by $\sigma_i$. It follows from Lemma 4q that the corollary is valid.**

**THEOREM G. Let G be a graph containing a hamilton circuit, $H_C$, while $H_0$ is a pseudo-hamilton circuit each of whose edges lies in $K_n$ - G. Then there always exists a set of successive $H_i$-admissible permutations, $\sigma_i$**

(i = 0, 1, 2, 3, ... ) to obtain $H_C$.
PROOF. Lemmas 2q and 3q prove the theorem when $H_0$ lies in $K_n$ - G.

THEOREM GEN. Let G be a graph containing a hamilton circuit, $H_C$, while $H_0$ is an *arbitrary* pseudo-hamilton circuit. Then we can always use successive $H_i$-admissible permutations to obtain $H_C$ *provided we may assume at most once that each arc vertex of* $H_0$ *is a pseudo-arc vertex* when we construct the $H_i$-admissible permutations.

PROOF. The proof is the same as that of Theorem G except that, in some cases, we must use arc vertices instead of pseudo-arc vertices *at most once* to eventually obtain $H_C$.

COMMENT. The limitation of Theorem G with regards to $H_0$ do not necessarily invalidate the use of Algorithm G. Even if $H_0$ contains some arc vertices, using the contracted graph, G', we can generally pass through the original arc vertices of $H_0$' as the algorithm proceeds. The best procedure would be to start by using the pseudo-arc vertices of $H_0$' in constructing admissible permutations. If necessary when constructing admissible permutations, we consider the arc vertices as pseudo-arc vertices the first time we use them. Using randomly chosen rotations in Algorithm G as well as arc vertices (if necessary) generally gives us the opportunity to go through each pseudo-arc vertex of $H_0$' in 2nlogn iterations of the algorithm. The second phase of the algorithm (to obtain a hamilton circuit in G') generally requires another 2nlogn iterations. We note that the existence Theorem G and Theorem Gen do not require the use of backtracking or rotations.

COROLLARY 1. Let $H_0$' be a pseudo-hamilton circuit of the contracted graph, G', of G, where all of the edges of $H_0$' lie in $K_n$ - G'. Then a necessary condition for G to contain a hamilton circuit is that each vertex, v, of G' lies on an $H_0$'-admissible permutation of G'.

PROOF. From Theorem G, if $H_0$' lies in $K_n$ - G', then no arcs of the hamilton circuit, $H_C$', lie on $H_0$'. Thus, we may use any vertex, v, of G' to start the algorithm given in Theorem G, i.e., v , lie on an $H_0$'-admissible permutation.

In Corollary 2, assume that we delete the use of backtracking and rotations from Algorithm G and Algorithm D. Call the algorithms obtained Algorithm $G_{simp}$ and Algorithm $D_{simp}$, respectively.

COROLLARY 2. Given a graph, G, or a directed graph, D, assume that we construct a pseudo-hamilton circuit, $H_0$, in G (D). Then, using quantum computing, if G (D) contains a hamilton circuit, we can obtain every hamilton circuit in G (D) using Algorithm $G_{simp}$ (Algorithm $D_{simp}$) where those branches allow us to proceed only by backtracking are **eliminated. In general, it requires no more than [.5n] iterations along a branch to obtain any circuit.**

**Note. If $H_0$ contains one or more arc vertices, we must assume that each arc vertex, v, *may* be considered a pseudo-arc vertex, the *first* time we choose an arc or pseudo-arc emanating from it. Once we have gone through each arc vertex of $H_0$ once or have assumed that particular arcs on $H_0$ belong to *all* hamilton circuits obtainable from a fixed branch, then we may assume that *only* pseudo-arcs may be chosen henceforth. This cuts down on the number of branches we have to consider.**

**Since we don't use rotations or backtracking in Corollary 2, PSEUDO consists only of vertices.**

**PROOF OF COROLLARY 2. When we start the proof, we assume (*whether or not this is really the case*) that all vertices of $H_0$ are pseudo-arc vertices. Thus, $PSEUDO_0 = V$. We now place the each pseudo-arc $[v, v_i]$ of $H_0 \cap G$ on a separate branch with i = 0, 1, 2, ... , n. Here $v_i \neq v_j$ if $i \neq j$. If $[v, v_i]$ lies on a branch, we can make no further choices of an arc out of v on that branch. Using all possible arcs in G in conjunction with each $[v, v_i]$, we test for $H_0$-admissibility. Each case of $H_0$-admissibility to obtain all possible pseudo-hamilton circuits, $H_{0\alpha_k}$ Here**

$$\alpha = 1, 2, \ldots \; ; k = 1$$

**As in Algorithms G and D, we use abbreviations to shorten computation. Each branch now has a separate set, $PSEUDO_{\alpha_k}$, associated with it. Any given $H_{0\alpha_k}$ contains at least two arc vertices. If k = 2, each $H_{0\alpha_k}$ contains at least four arc vertices, and we have deleted at least four vertices from $PSEUDO_{\alpha_k}$. We eliminate all branches that require backtracking to continue. It follows that after at most [.5n] iterations, we obtain all hamilton circuits in G.**

**THEOREM H. Let G be a random graph of degree n which contains a hamilton circuit. Then if Algorithm G is allowed to run through $12n^3 \log n$ iterations, the probability of obtaining a hamilton circuit in G is at least**

$$\left(1 - \frac{1}{n^{3n^2}}\right)^2$$

**The running time of the algorithm is at most MG where**

$$\begin{aligned} MG = {} & (192(\log n)^2 + 24\log n + 24)[n^{3.5}] \\ & + (336(\log n)^3 + 576(\log n)^2 + 156\log n + 4)[n^3] \\ & + (3\log n + 1)n^2 + (5\log n)n + \log n \end{aligned}$$

**PROOF. Given a suitably large number of iterations, by the Law of Large Numbers, the number of successes approaches the number of iterations multiplied by the probability of success yielding**

$$\left(\frac{13m^4 - 8m^3 - 6m^2 - 1}{16m^4}\right)(6n^3 \log n) > 3n^3 \log n$$

The probability of going through a fixed vertex is $\frac{1}{n-1}$. It follows that the probability of going through each vertex using $6n^3\log n$ iterations is at least

$$1 - \left(1 - \frac{1}{n}\right)^{3n^3 \log n} > 1 - \frac{1}{n^{3n^2}}$$

The probability that each edge obtained from a pseudo-3-cycle lies in G is at least $\frac{1}{n}$. Therefore, the probability that by using another $6n^3\log n$, we can obtain a hamilton cycle from a hamilton path is at least

$$\left(1 - \frac{1}{n^{3n^2}}\right)^2$$

We now obtain the running time of the algorithm.

(I) Constructing G. It takes nlogn r. t. to construct a balanced, binary search tree, $G_{TREE}$, with n branches. It requires at most $(\log n)^2$ r. t. to randomly pick a vertex from a V X V adjacency matrix of $K_n$ and logn r. t. to place a vertex on a branch. Since the number of edges in the graph is at most $\binom{n}{2}$, it takes at most $\underline{n^2 \log n}$ r. t. to place the vertices of each edge on the branches it is incident to - on $G_{TREE}$. It takes at most $\underline{2\binom{n}{2}}$ r. t. to use a counting number on each branch to count the number edges incident to the vertex on that branch.

TOTAL R.T. : $\underline{(\log n + 1)n^2 - n}$

(II) Searching for vertices of degree 2. It takes nlogn r. t. to obtain each vertex of degree 2 and place it in a balanced, binary search tree, TWO. We use the following procedure to construct each 2-path, $S_\alpha$. Using TWO, let b be an arbitrary vertex of degree 2 with unique edges [a, b] and [b, c] incident to it. If deg(c) > 2, stop. Otherwise, obtain the edge of G incident to c other than [b, c], say [c, d]. If deg(d) > 2, stop, etc.. Now we consider a. If deg(a) > 2, stop. Otherwise, we consider the edge of G other than [a, b] which is incident to a, say [z, a]. If deg(z) > 2, stop. Otherwise, given the edge [y, z], consider deg(y), etc. It requires at most 6logn r.t. to obtain

$$S_\alpha = \{y, z, a, b, c, d\}$$

**Generally, it requires no more than (r+1)(logn) r. t. to find and place each $S_\alpha$ of length r in the balanced, binary search tree, TWOPATH. Even if the entire graph consists only of vertices of degree 2, it requires at most <u>(n+1)(logn)</u> to obtain each 2-path and place it in TWOPATH.**
**TOTAL R. T.: <u>(logn)n + logn</u>**

**(III) Replacing each subpath by a 2-vertex to form G'. Using $G_{TREE}$, we construct a duplicate of it, $G_{TREE}'$ leaving out branches headed by vertices of degree 2 obtained from TWO. We then obtain each 2-path from TWOPATH and construct a 2-vertex, say vv*, corresponding to its endpoints. We then obtain from $G_{TREE}$ the edges incident to either v or v* and place them on two distinct vv* branches. This yields G'. The number of degrees of edges in G' is no greater than n(n-1). The maximum number of vertices of degree 2 is n. It takes at most nlogn r. t. to delete each vertex of degree 2 and the edges incident to it from $G_{TREE}'$. On the other hand, it requires at most (n + 1)logn r. t. to place each 2-vertex on $G_{TREE}'$ and n(n-1)logn r. t. to add two branches of vertices incident to each vertex to obtain G'. Thus, the total r. t. for (III) is less than <u>2n(n-1)logn</u>. For simplicity, when placing a 2-vertex on G', its order of magnitude is always that of the smaller natural number.**
**TOTAL R. T. : <u>(2logn)$n^2$ + (-2logn)n</u>**

**(IV) Given G', construct $H_0'$, ORD($H_0'$). It requires at most nlogn r. t. to construct a tree containing the natural numbers from 1 through n, say N. It requires logn r. t. to do each of the following: pick a random number from N, place the number in the balanced, binary search tree, $H_0'$, delete the number from N, rebalance the tree. It requires logn r. t. to place the ordinal value of the number, ORD(1), in ORD($H_0'$) and rebalance the tree. Thus, it requires at most 6nlogn r. t. to construct both $H_0'$ and ORD($H_0'$).**
**TOTAL R. T.: <u>(6logn)n</u>**

**(V) We now compute the running time of [$n^{.5}$] iterations using Algorithm G when n is a fixed natural number and $H_\alpha'$ is a pseudo-hamilton circuit with α of the form i[$n^{.5}$]. We note that $PSEUDO_\alpha$ contains the actual values of vertices – not their ordinal values with respect to a fixed pseudo-hamilton circuit.**

**(1) Choose two pseudo-arcs from $PSEUDO_{\alpha+j}$, say [a, r], [b, s].**

**If n = $e^r$, the number of bits in the four numbers is no greater than <u>4logn</u>.**

**(2) Find the ordinal numbers in ORD($H_\alpha'$) of a and b.**
**This requires at most <u>2logn</u> r.t. .**

**(3) Locate ORD(a) on $A_{\alpha+j}$ and determine whether or not "-" precedes it. Do the same thing for ORD(b).**
**This requires at most <u>2j</u> r.t. .**

**(4) W.l.o.g., assume that neither ordinal value is preceded by "-". Thus, we can assume that a and b are both pseudo-arc vertices.**

If one or more of them had a "-" in front of it, it would take at most 2J r. t. to obtain both ORD(r) and ORD(s).

(5) Choose arcs respectively out of ORD(a) and ORD(b), say [a, p], [b, q].
This requires at most 2logn r.t. .

(6) Obtain ORD(p) and ORD(q).
This requires at most 2logn r.t. .

(7) Determine the order of ORD(a), ORD(b), ORD(p) and ORD(q) along $A_{\alpha}$.
This requires at most 4j r.t. .

(8) If the numbers interlace, then
[ORD(a), $H_{\alpha+j}$' (ORD(p))], [ORD(b), $H_{\alpha+j}$' (ORD(q))]
define an $H_{\alpha+j}$ '-admissible permutation, $s_{\alpha+j}$. Given that the ordinal numbers interlace, we now construct $A_{\alpha+j+1}$.
This requires at most 8j r. t. which occurs when $s_{\alpha+j}$ is a POTDTC. Otherwise, it requires at most 6j r. t. .

(9) We test for at most $2(\log n)^2$ possible admissible permutations.
This requires at most $28(\log n)^3 + 32j(\log n)^2$ r. t.

(10) Given that all cases yield admissible permutations, it requires at most $(8)2(\log n)^2 = 16(\log n)^2$ r. t. to obtain all cases in which

M = (no. of pseudo-arc vertices) – (no. of pseudo-arcs)

in an admissible permutation is greater than 0 and to obtain the largest value of M.

(11) Given that all cases in (10) have the same value for M, it requires at most $16(\log n)^2$ r. t. to obtain a permutation from among these, one in which the sum of the degrees of its vertices is greatest.

(12) Thus, the r. t. up to this point in an iteration of Algorithm G is at most

$$28(\log n)^3 + 32(\log n)^2 (j + 1)$$

(13) It requires at most logn r. t. to randomly choose an arc out of a vertex, say v, in order to perform a rotation.

(14) It requires at most 4j r. t. to perform a rotation.

(15) It requires at most 8logn r. t. to delete four pseudo-arcs from $PSEUDO_{\alpha+j}$. It requires at most 4logn r. t. to place two pseudo-arcs in $PSEUDO_{\alpha+j}$ thus forming PSEUDO $_{\alpha+j+1}$.

(16) It thus takes at most

$$28(\log n)^3 + (32(\log n)^2 + 4)j + 32(\log n)^2 + 13\log n$$

r. t. up to this point in (V).

Rewriting the above,

$$\sum_{j=1}^{j=[n^{.5}]} \{(32(\log n)^2 + 4)(j) + (28(\log n)^3 + 32(\log n)^2 + 13\log n)$$

(17) We now use $A_{\alpha+[n^{.5}]}$ and ORD($H_\alpha$') to construct $H_{\alpha+[n^{.5}]}$' and ORD($H_{\alpha+[n^{.5}]}$').

It takes logn r. t. to find the vertex corresponding to each ordinal value of $A_{\alpha+[n^{.5}]}$ using $H_\alpha$'. Once we have obtained each ordinal value, we start the construction of $H_{\alpha+[n^{.5}]}$' and ORD($H_{\alpha+[n^{.5}]}$'). This requires at most 2nlogn r. t. . Therefore,

TOTAL R. T. for (V) is

$$\sum_{j=1}^{[n^{.45}]} \{(32(\log n)^2 + 4)j + (28(\log n)^3 + 32(\log n)^2 + 13\log n)\} + 2n\log n$$

$$= (16(\log n)^2 + 2\log n + 2)n + (28(\log n)^3 + 48(\log n)^2 + 13\log n + 2)[n^{.5}]$$

Since we will apply $12n^{2.5}$ cases of $[n^{.5}]$ iterations,

TOTAL R. T. for this phase of the algorithm is

$(192\log^2 n + 24\log n + 24)n^{3.5}$

$+ (336\log^3 n + 576\log^2 n + 156\log n + 4)n^3$

(VI). Replacing $H_C$' by $H_C$. Given a 2-vertex, vv*, we check $A_C$, the last abbreviation constructed, to see if there is a "-" in front of vv*. This requires at most $[n^{.5}]$ r. t. . If there is such a sign, we must go to the last number in the subpath and add it to $H_C$' in place of (w.l.o.g.) v*. We then must obtain the next to last number and add it to $H_C$' after v*. We continue in this manner until we reach v. If the length of the subpath is r, it requires at most (r+1)logn r. t. to replace vv* by the subpath it represent. It also requires at most $\binom{r}{2}$ running time to obtain to obtain each vertex. Thus, the maximum r. t. necessary is never greater than (n+1)logn + $\binom{n}{2}$ which simplifies to $n^2 + (\log n - .5)n + \log n$. Thus,

TOTAL R. T. : $n^2 + (\log n - .5)n + \log n$

Thus, the running time, MG, of the algorithm using a total of $12n^3$ iterations is at most

$(192(\log n)^2 + 24\log n + 24)[n^{3.5}]$

$$+ (336(\log n)^3 + 576(\log n)^2 + 156\log n + 4)[n^3]$$
$$+ (3\log n + 1)n^2 + (5\log n)n + \log n$$

**THEOREM I. Let D be a random directed graph of degree n which contains a hamilton cycle. Then Algorithm D obtains a hamilton cycle in D with probability at least**

$$(1 - \frac{1}{n^{3n^2}})^2$$

**Its running time is at most MG.**

**PROOF. The algorithm is essentially the same as in Theorem H except that we don't use rotations but we do use backtracking. The probability for a successful iteration is**

$$\frac{13m^4 - 8m^3 - 6m^2 - 1}{16m^4}$$

**while the probability for failure is**

$$\frac{3m^2 + 8m^3 + 6m^2 + 1}{16m^4}$$

**Thus, the net probability for success, $p_{net}$, after backtracking is taken into account is**

$$\frac{10m^4 - 16m^3 - 12m^2 - 2}{16m^4}$$

**Therefore,**

$$(p_{net})(6n^3\log n) > 3n^3\log n$$

**That is, for a large number of iterations, the number of successful cases is greater than $3n^3$logn, the same number as that used in the previous case. Furthermore, the number of iterations (including those requiring backtracking) is still $12n^3$logn. Furthermore, we don't require the construction of a rotation. It follows that the running time for the algorithm is no greater than MG.**

**The algorithm doesn't require constructing rotations. On the other hand, it does require backtracking.**

**Before going on to Conjecture A, we prove the following theorem.**

**THEOREM R. Let G' be the contracted graph of an arbitrary graph G where G contains a hamilton circuit. Suppose H' is a pseudo-hamilton circuit of G'. Assume that the following are true:**

**(1) The vertices of H' are equally spaced along a circle C.**

**(2) a is a pseudo-arc vertex of H'.**

**(3) $e_1$ = [H'(a), H'(d)] and $e_2$ = [d, H'(e)] are edges of G'- H' which do not intersect in C.**

**(4) [a, x] is an edge of G' - H' determining a rotation, r, of H' such that H'(a) is a pseudo-arc vertex of H* = H'r. Then if – going in a clockwise direction – x lies between H'(e) and d, $e_1$ and $e_2$ intersect in H* and determine an H*-admissible pseudo-3-cycle.**

**PROOF. W.l.o.g., let**
**H' = (1 2 3 ... *21* ... *30* ... *39* *40* 41 ... n).**
**Assume that a = 20, H'(a) = 21, H'(d) = 40, d = 39, H'(e) = 30,**
**x = 35. Then H* = (1 2 ... 20 35 34 33 32 31 *30*...**
**_21_ 36 37 38 _39_ _40_ ... n).**
**We note that r *reversed* the order of H'(a) = 21 and H'(e) = 30. Thus, 30, 21, 39, 40 interlace the vertices of $e_1$ and $e_2$, i.e., $e_1$ and $e_2$ intersect in the circle defined by H*. It follows that since H'(a) is a pseudo-arc vertex by hypothesis, the latter two edges define a H*-admissible pseudo-3-cycle.**

**COROLLARY. Let $e_1$ and $e_2$ be defined as in Theorem R. Define $e_3$ = [H'(H'(d)), g] = [f, g]. Assume that r = [a, x] defines a rotation. Then if**

**(a) if x < d and $e_1$ and $e_2$ intersect,**

**or**

**(b) if x > f and $e_1$ and $e_3$ intersect,**

**they define an H*-admissible pseudo-3-cycle.**

**PROOF. [a] On any circle H* in which --- going in a clockwise direction --- d is the predecessor of H'(d) and $e_1$ and $e_2$ intersect, they define an H*-admissible pseudo-3-cycle.**

**[b] If x > f and $e_1$ and $e_3$ intersect, r changes f from a successor of H'(d) going clockwise along H' to a predecessor of H'(d) going clockwise along H*. Thus, the two edges define an H*-admissible pseudo-3-cycle.**

**CONJECTURE. Let G be an arbitrary graph which contains a hamilton circuit. Let G' be its contracted graph. Then, using Algorithm G, we can always obtain a hamilton circuit of G in polynomial time.**

**Comment. The reasoning here is that the random choice of edges in the rotation at the end of each iteration gives us the means to go successfully through each vertex of V in polynomial time even if we don't obtain an admissible permutation in each iteration. Theorem R and its corollary show that it is possible by using a rotation to convert two edges which don't define an admissible permutation into a pair defining an admissible permutation. Finally, we note that in Theorems H and I, the probability of failure decreases exponentially as the running time increases polynomially.**

In general, it is useful to cut down on the number of arc vertices in $H_0$. To do this, we first define a *$cH_0$-admissible permutation.* An $H_0$-admissible permutation the majority of whose edges lie in $K_n$ - G' is called a $cH_0$-admissible permutation.

The restriction on the number of arc vertices allowed in $H_0$ ' is not that difficult to deal with *provided that the number of edges in G' is considerably smaller than the number in its complement, $K_n$ - G'*. We fist note that Theorem G only requires that G' (and therefore G) contain a hamilton circuit. Thus, we don't have to randomly construct $H_0$ '. Any good heuristic or the algorithm given in this paper may be used to construct $H_0$ '. If we wish to obtain $cH_i$' = $H_0$ ' with a minimum number of arcs in G', we could use quantum computing to eliminate edges in G' from $cH_0$ ', replacing them with pseudo-arcs of G' until we obtain an $H_0$ ' containing close to the minimum possible number of arcs vertices in any possible $H_0$ '.

EXAMPLE. Let G be the set of all edges, [a, b], in $K_{32}$ such that one vertex is even and the other is odd. Assume that all vertices *not* in italics and underlined are *arc vertices* of G. Define $cPSEUDO_i$ as the set of arc vertices of G lying in $H_i$.

Let

$H_C$ = (1 2 3 4 5 6 7 8 9 10 11 12 13 14 15 16 17 18 19 20
21 22 23 24 25 26 27 28 29 30 31 32)

be a hamilton circuit in G. Define

$cH_0$ = (*1* 22 *28* 29 *21* 12 *8* 15 25 27 *19* 6 *4* 23 *7* *24*
*31* *18* 11 13 *5* *2* *9* 32 30 20 *26* *17* 14 *16* *3* *10*)

$cPSEUDO_0$ = (1, 28, 21, 8, 19, 4, 7, 24, 31, 18, 5, 2, 9, 26, 17, 16, 3, 10}

If $s_0$ = (1 28 6) then

$cH_1$ = (1 29 *21* 12 *8* 15 25 27 *19* *6* 22 28 *4* 23 *7* *24* *31* *18* 11 13
*5* *2* *9* 32 30 20 *26* *17* 14 *16* *3* *10*)

$cPSEUDO_1$ = (21, 8, 19, 4, 7, 24, 31, 18, 5, 2, 9, 26, 17, 16, 3, 10}}

If $s_1$ = (21 24) (4 3) then

$cH_2$ = (1 29 21 *31* *18* 11 13 *5* *2* *9* 32 30 20 *26* *17* 14 *16* 3 23 *7*
24 12 *8* 15 25 27 *19* 6 22 28 4 *10*)

$cPSEUDO_2$ = {31, 18, 5, 2, 9, 26, 17, 16, 7, 8, 19, 10}

If $s_2$ = (5 26)(9 8), then

$cH_3$ = (1 29 21 *31* *18* 11 13 5 *17* 14 *16* 3 23 *7* 24 12 8 32 30 20 26
*2* 9 15 25 27 *19* 6 22 28 4 *10*)

$cPSEUDO_3$ = {31, 18, 17, 16, 7, 2, 19, 10}

If $s_3$ = (31 2 4), then

$cH_4$ = (1 29 21 31 9 15 25 27 *19* 6 22 28 4 *18* 11 13 5 *17* 14 *16* 3 23 *7*
24 12 8 32 30 20 26 2 *10*)

$cPSEUDO_4$ = {19, 18, 17, 16, 7, 10}

**If $s_4$ = (18 14 10), then**

**$cH_5$ = (1 29 21 31 9 15 25 27 *19* 6 22 28 4 18 *16* 3 23 7 11 13 5 *17* 14 24 12 8 32 30 20 26 2 *10*)**

**$cPSEUDO_5$ = {19, 16, 17, 10}**

**If $s_5$ = (19 6 2), then**

**$cH_6$ = (1 29 21 31 9 15 25 27 19 3 23 7 11 13 5 *17* 14 24 12 8 32 30 20 26 2 6 22 28 4 18 16 *10*)**

**$cPSEUDO_6$ = {17, 10}**

**$cH_6$ can be used as $H_0$. We cannot obtain a pseudo-hamilton circuit with fewer than two pseudo-arcs: There always must be a change from an even number to an odd number and from an odd number to an even number.**

## VI. A HEURISTIC FOR THE SYMMETRIC TRAVELING SALESMAN PROBLEM.

**In this section, we discuss a heuristic for the symmetric traveling salesman problem using Algorithm G.**

**(1) Construct an n X n symmetric matrix. say N, containing only 0 entries along its diagonal.**
**(2) Randomly construct an n-cycle, say $h_0'$. Let the corresponding pseudo-hamilton circuit be $H_0'$.**
**(3) Assign a positive weight, w(i,j), to each non-diagonal entry, (i,j).**
**(4) Apply a simple heuristic (say, Lin-Kernigan [18]) to $H_0'$ to obtain an approximation to a smallest sum of weights of an n-cycle, say $H_0$.**

**We now make the following definitions:**
**The *defining arcs of an admissible permutation* are two arcs associated with the permutation which intersect. A set of arcs is *good* if the sum of their weights is less than the sum of the weights of arcs on $H_i$ which have the same respective initial vertices. Then the following hold: *If the sum of the weights of the arcs associated with an $H_i$-admissible permutation is less than the sum of the weights of the arcs with corresponding initial vertices, then the permutation is good. Furthermore, if a rotation determined by the arc [x,y] has the property that***

$$w[x,y] + w[H_i(x), H_i(y)] < w[x, H_i(x)] + w[y, H_i(y)]$$

***then the rotation is good.***

**(5) We next make rotations through each of the vertices of $H_i$ until we can no longer make a good rotation.**
**(6) An iteration now consists of one of the following:**
   **(a) a good admissible permutation which is *not* followed by a rotation.**
   **(b) for a 3-cycle, an admissible permutation such that the sum of the weights defined by a set of defining arcs which added to the sum of the weights of a rotation define a *good* set of arcs.**

(c) for a POTDT, an admissible permutation for which the following is true:

(i) the sum of the weights of a set of defining arcs, added to the sum of the weights of a rotation through the third arc, added to the sum of the weights of a rotation through the fourth arc, define a *good* set of arcs;

(ii) the sum of the weights of a defining set of arcs and a third arc, added to the sum of the weights of a rotation through the fourth vertex, define a *good* set of arcs.

(d) If $H_i$ is a weighted n-cycle before 6 (a), (b) and (c) have been applied to it, while $H_{i+1}$ is the result after such applications, then

(I) in the case of a 3-cycle,

$$W(H_i) - W([a, H_i(a)]) > W(H_{i+1}) - W([H_i(a), H_{i+1}(H_i(a))]);$$

(II) in the case of a POTDTC,

$$W(H_i) - \{W([c, H_i(a)]) + W([d, H_i(b)])\}$$
$$> W(H_{i+1}) - \{W([H_i(a) - H_{i+1}(H_i(a))]) + W([H_i(b) - H_{i+1}(H_i(b))]\}$$

Note. We require at least O(i) bits to test a sequence of an admissible permutation followed by one or two rotations on $A_i$ to see if a good set of arcs is defined. We then must delete arcs after the testing is done.

(7) As we proceed in (6), if $W(H_{i_1}) < W(H_0)$, we place $W(H_{i_1})$ and $H_{i_1}$ in a queue. If $W(H_{i_2}) < W(H_{i_1})$, we delete $W(H_{i_1})$ and $H_{i_1}$ from the queue and replace them with $W(H_{i_2})$ and $H_{i_2}$. The algorithm concludes when there is no vertex out of which we can construct a good $H_{i_\alpha}$-admissible permutation according to the rules given in (6) with $W(H_{i_\alpha}) > W(H_{i_{\alpha+1}})$.

## REFERENCES.

[1] M. Ajtai, J. Kömlos, E. Szeremédi, First occurrence of hamilton cycles in random graphs, Annals of Discrete Mathematics 27 (1985), 173-178.

[2] D. Angluin and L. Valiant, Fast probabilistic algorithms for Hamilton circuits and matchings, J. Comp. System Sci. (1978), 155-193.

[3] B. Bollobás, The evolution of sparse graphs in "Graph Theory and Combinatorics Proceedings. Cambridge Combinatorial Conference in Honor of Paul Erdös, 1984" (B. Bollobás, Ed.) pp. 335-357.

[4] B. Bollobás, T. I. Fenner and A. M. Frieze, An algorithm for finding hamilton paths and cycles in random graphs,"Proceedings, 17th Annual ACM Symposium, 1985," pp. 430-439.

[5] B. Bollobás, *Random Graphs*, Academic Press, New York, London, 1985, p. 152.

[6] Cooper, C. and Frieze, A. M., Hamilton cycles in a class of random directed graphs, J. Combin. Theory, Ser. B. 62 (1994), no. 1, 157-163.

[7] Dickson, L. E., History of the Theory of Numbers, Vol. 2, Chelsea Publishing Company, New York (1966), 4-5.

**[8] Erdös, P. and Renyi, A., On the strength of connectedness of a random graph, Acta Math. Acad, Sci, Hung. 12 (1961), 261-267.**
**[9] Feller, W., *An Introduction to Probability Theory and its Applications, Vol. 1,* Wiley Publications in Statistics, 2nd Edition, John Wiley & Sons, New York, London, 1950, pp.160-161.**
**[10] T. I. Fenner and A. M. Frieze, On the existence of hamilton cycles in a class of random graphs, Discrete Math. 45, 55-63**
**[11] T. I. Fenner and A. M. Frieze, On the connectivity of random m-orientable graphs and digraphs, Combinatorica 2 (4) (1982), 347-359.**
**[12] A. Frieze, An algorithm for finding hamilton cycles in random directed graphs, Journal of Algorithms 9 (1988), 181-204.**
**[13] Frieze, A. M., Luczak, Tomasz, Hamiltonian cycles in a class of directed graphs: one step further. Random solution for the existence of Hamilton cycles in random graphs, Discrete Math. 43 (1983), 55-63.**
**[14] A. Frieze, M. Jerrum, M. Molloy, R. Robinson, N. Wormald, Generating and counting hamilton cycles in random regular graphs, Journal of Algorithms, Vol. 21 (1996), no. 1, 176-198.**
**[15] W. Hoeffding, Probability inequalities for sums of bounded, random variables, J. American Statistical Assoc. (1963), 13-30.**
**[16] D. Knuth, *The Art of Computer Programming, Vol. 3,* Addison-Wesley, Reading, Mass., 1973.**
**[17] J. Komlös and E. Szemerédi, Limit distribution for the existence of Hamilton cycles in random graphs, Discrete Math. 43 (1983), 55-63.**
**[18] S. Lin, B. W. Kernighan, An effective, heuristic method for the traveling salesman problem, Operations Res. 21 (1973), 498-517.**
**[19] C. J. H. McDiarmid, Clutter percolation and random graphs, Math. Programming Studies, vol. 13 (1980), pp. 17-25.**
**[20] Palásti, I., On the strong-connectedness of directed random graphs, Math. Programming Stud., vol. 13 (1980), pp. 17-25.**
**[21] Nicholas C. Wormald, The asymptotic connectivity of random regular graphs, J. of. Comb. Theory, Series B, Vol. 31, (1981), 156-167.**

**Current Address: 188-83 85th Road**
**Holliswood, New York 11423**
**e-mail addresses: HowardKleiman@qcc.cuny.edu**
**hkleiman@banet.net**